\documentclass [12pt]{article}
\usepackage{amsmath,amssymb}

\begin{document}

\title{Explicit representation formulas for the minimum norm least squares solutions of some  quaternion
matrix equations.}
\author{Ivan Kyrchei \footnote{Pidstrygach Institute for Applied Problems of Mechanics and Mathematics of NAS of Ukraine,
str. Naukova 3b, Lviv, 79060, Ukraine, kyrchei@lms.lviv.ua}}
\date{}
 \maketitle

\begin{abstract}
Within the framework of the theory of the column and row
determinants, we obtain explicit representation formulas (analogs
of Cramer's rule) for the minimum norm least squares solutions of
quaternion matrix equations ${\bf A} {\bf X} = {\bf B}$, $ {\bf X}
{\bf A} =  {\bf B}$ and ${\bf A} {\bf X} {\bf B} = {\bf D} $.
\end{abstract}



\section{Introduction}
\newtheorem{theorem}{Theorem}[section]
\newtheorem{lemma}{Lemma}[section]
\newtheorem{definition}{Definition}[section]
\newtheorem{proposition}{Proposition}[section]
\newtheorem{remark}{Remark}[section]
In recent years quaternion matrix equations have been investigated
by many authors (see, e.g., \cite{ji1}-\cite{ky2}). For example,
 Jiang, Liu, and Wei \cite{ji1} studied the solutions of the
general quaternion matrix equation ${\rm {\bf A}}{\rm {\bf X}}{\rm
{\bf B}} - {\rm {\bf C}}{\rm {\bf Y}}{\rm {\bf D}}={\rm {\bf E}}$,
and Liu \cite{li} studied the least squares Hermitian solution of
the quaternion matrix equation $({\rm {\bf A}}^{H}{\rm {\bf
X}}{\rm {\bf A}},{\rm {\bf B}}^{H}{\rm {\bf X}}{\rm {\bf
B}})=({\rm {\bf C}}{\rm {\bf D}})$, Wang, Chang, and Ning
\cite{wa5} derived the common solution to six quaternion matrix
equations, etc.

However the understanding of the problem for determinantal
representing the least squares  solutions of the quaternion matrix
equations, in particular
\begin{equation}\label{eq:def_AX=B}
  {\rm {\bf A}}{\rm {\bf X}}={\rm {\bf
B}},
\end{equation}
 \begin{equation}\label{eq:def_XA=B}
 {\rm {\bf X}}{\rm {\bf A}}={\rm {\bf B}},
\end{equation}
 \begin{equation}\label{eq:def_AXB=D}
 {\rm {\bf A}}{\rm
{\bf X}}{\rm {\bf B}} ={\rm {\bf D}},
\end{equation}
 has not yet reached a
satisfactory level. The reason was the lack of appropriate
noncommutative determinant. Many authors had tried to give the
definitions of the determinants of a quaternion matrix, (see, e.g.
\cite{as}-\cite{ge2}). Unfortunately, by their definitions it is
impossible for us to give an determinant representation of an
inverse matrix.

But in \cite{ky1}, we defined the row and column determinants and
the double determinant  of a square matrix over the quaternion
skew field. As applications we obtained the determinantal
representations of an inverse matrix by an analogue of the adjoint
matrix in \cite{ky1} and the Moore-Penrose inverse  over
quaternion skew field ${\mathbb{H}}$ in \cite{ky3}. In \cite{ky2}
we also gave  Cramer's rule for the solution of nonsingular
quaternion matrix equations (\ref{eq:def_AX=B}),
(\ref{eq:def_XA=B}) and (\ref{eq:def_AXB=D}) within the framework
of the theory of the column and row determinants. In
\cite{so1,so2}, the authors obtained the generalized Cramer rules
for the unique solution of the matrix equation
(\ref{eq:def_AXB=D}) in some restricted conditions within the
framework of the theory of the column and row determinants as
well.

In this paper we  aim to obtain explicit representation formulas
(analogs of Cramer's rule) for the minimum norm least squares
solutions of quaternion matrix equations (\ref{eq:def_AX=B}),
(\ref{eq:def_XA=B}) and (\ref{eq:def_AXB=D}) without any
restriction.  The paper is organized as follows. In Section 2, we
start with some basic concepts and results  from the theory of the
column and row determinants which are necessary for the following.
The theory of the column and row determinants of a quaternionic
matrix is considered completely in \cite{ky1}. In Section 3, we
give the theorem about determinantal representations of the
Moore-Penrose inverse over the quaternion skew field derived in
\cite{ky3}. In Section 4, we obtain explicit representation
formulas  for the minimum norm least squares solutions of
quaternion matrix equations (\ref{eq:def_AX=B}),
(\ref{eq:def_XA=B}) and (\ref{eq:def_AXB=D}). In Section 5, we
show a numerical example to illustrate the main result.

\section{Elements of the theory of the column and row determinants.}
Throughout the paper, we denote the real number field by ${\rm
{\mathbb{R}}}$, the set of all $m\times n$ matrices over the
quaternion algebra
\[{\rm {\mathbb{H}}}=\{a_{0}+a_{1}i+a_{2}j+a_{3}k\,
|\,i^{2}=j^{2}=k^{2}=-1,\, a_{0}, a_{1}, a_{2}, a_{3}\in{\rm
{\mathbb{R}}}\}\] by ${\rm {\mathbb{H}}}^{m\times n}$  and   its
subset of matrices of rank $r$ by ${\rm {\mathbb{H}}}_{r}^{m\times
n} $.  Let ${\rm M}\left( {n,{\rm {\mathbb{H}}}} \right)$ be the
ring of $n\times n$ quaternion matrices.

The \textit{conjugate} of a quaternion $a = a_{0} + a_{1} i +
a_{2} j + a_{3} k \in {\rm {\mathbb{H}}}$ is defined by $\overline
{a} = a_{o}-a_{1}i- a_{2}j - a_{3}k$. The \textit{Hermitian
adjoint matrix} of ${\rm {\bf A}} = \left( {a_{ij}} \right) \in
{\rm {\mathbb{H}}}^{n\times m}$ is called the matrix ${\rm {\bf
A}}^{ *} = \left( {a_{ij}^{ *} } \right)_{m\times n} $ if
$a_{ij}^{ *} = \overline {a_{ji}}  $ for all $i = 1,\ldots,n $ and
$j = 1,\ldots,m$.
 The matrix ${\rm {\bf A}} = \left( {a_{ij}}  \right) \in {\rm
{\mathbb{H}}}^{n\times m}$ is \textit{Hermitian} if ${\rm {\bf
A}}^{ *}  = {\rm {\bf A}}$.

Suppose $S_{n}$ is the symmetric group on the set
$I_{n}=\{1,\ldots,n\}$.
\begin{definition}
 The $i$th row determinant of ${\rm {\bf A}}=(a_{ij}) \in {\rm
M}\left( {n,{\mathbb{H}}} \right)$ is defined by
 \[{\rm{rdet}}_{ i} {\rm {\bf A}} =
{\sum\limits_{\sigma \in S_{n}} {\left( { - 1} \right)^{n -
r}{a_{i{\kern 1pt} i_{k_{1}}} } {a_{i_{k_{1}}   i_{k_{1} + 1}}}
\ldots } } {a_{i_{k_{1} + l_{1}}
 i}}  \ldots  {a_{i_{k_{r}}  i_{k_{r} + 1}}}
\ldots  {a_{i_{k_{r} + l_{r}}  i_{k_{r}} }}\] \noindent for all $i
= 1,\ldots,n $. The elements of the permutation $\sigma$ are
indices of each monomial. The left-ordered cycle notation of the
permutation $\sigma$ is written as follows,
\[\sigma = \left(
{i\,i_{k_{1}}  i_{k_{1} + 1} \ldots i_{k_{1} + l_{1}} }
\right)\left( {i_{k_{2}}  i_{k_{2} + 1} \ldots i_{k_{2} + l_{2}} }
\right)\ldots \left( {i_{k_{r}}  i_{k_{r} + 1} \ldots i_{k_{r} +
l_{r}} } \right).\] \noindent The index $i$ opens the first cycle
from the left  and other cycles satisfy the following conditions,
$i_{k_{2}} < i_{k_{3}}  < \ldots < i_{k_{r}}$ and $i_{k_{t}}  <
i_{k_{t} + s} $ for all $t = 2,\ldots,r $ and $s =1,\ldots,l_{t}
$.
\end{definition}

\begin{definition}
The $j$th column determinant
 of ${\rm {\bf
A}}=(a_{ij}) \in {\rm M}\left( {n,{\mathbb{H}}} \right)$ is
defined by
 \[{\rm{cdet}} _{{j}}\, {\rm {\bf A}} =
{{\sum\limits_{\tau \in S_{n}} {\left( { - 1} \right)^{n -
r}a_{j_{k_{r}} j_{k_{r} + l_{r}} } \ldots a_{j_{k_{r} + 1}
i_{k_{r}} }  \ldots } }a_{j\, j_{k_{1} + l_{1}} }  \ldots  a_{
j_{k_{1} + 1} j_{k_{1}} }a_{j_{k_{1}} j}}\] \noindent for all $j
=1,\ldots,n $. The right-ordered cycle notation of the permutation
$\tau \in S_{n}$ is written as follows,
 \[\tau =
\left( {j_{k_{r} + l_{r}}  \ldots j_{k_{r} + 1} j_{k_{r}} }
\right)\ldots \left( {j_{k_{2} + l_{2}}  \ldots j_{k_{2} + 1}
j_{k_{2}} } \right){\kern 1pt} \left( {j_{k_{1} + l_{1}}  \ldots
j_{k_{1} + 1} j_{k_{1} } j} \right).\] \noindent The index $j$
opens  the first cycle from the right  and other cycles satisfy
the following conditions, $j_{k_{2}}  < j_{k_{3}}  < \ldots <
j_{k_{r}} $ and $j_{k_{t}}  < j_{k_{t} + s} $ for all $t =
2,\ldots,r $ and $s = 1,\ldots,l_{t}  $.
\end{definition}
The following lemmas enable us to expand ${\rm{rdet}}_{{i}}\, {\rm
{\bf A}}$ by cofactors
  along  the $i$th row and ${\rm{cdet}} _{j} {\rm {\bf A}}$
 along  the $j$th column respectively for all $i, j = 1,\ldots,n$.

\begin{lemma}\label{lemma:R_ij} \cite{ky1}
Let $R_{i{\kern 1pt} j}$ be the right $ij$-th cofactor of ${\rm
{\bf A}}\in {\rm M}\left( {n, {\mathbb{H}}} \right)$, that is, $
{\rm{rdet}}_{{i}}\, {\rm {\bf A}} = {\sum\limits_{j = 1}^{n}
{{a_{i{\kern 1pt} j} \cdot R_{i{\kern 1pt} j} } }} $ for all $i =
1,\ldots,n$.  Then
\[
 R_{i{\kern 1pt} j} = {\left\{ {{\begin{array}{*{20}c}
  - {\rm{rdet}}_{{j}}\, {\rm {\bf A}}_{{.{\kern 1pt} j}}^{{i{\kern 1pt} i}} \left( {{\rm
{\bf a}}_{{.{\kern 1pt} {\kern 1pt} i}}}  \right),& {i \ne j},
\hfill \\
 {\rm{rdet}} _{{k}}\, {\rm {\bf A}}^{{i{\kern 1pt} i}},&{i = j},
\hfill \\
\end{array}} } \right.}
\]
\noindent where  ${\rm {\bf A}}_{.{\kern 1pt} j}^{i{\kern 1pt} i}
\left( {{\rm {\bf a}}_{.{\kern 1pt} {\kern 1pt} i}}  \right)$ is
obtained from ${\rm {\bf A}}$   by replacing the $j$th column with
the $i$th column, and then by deleting both the $i$th row and
column, $k = \min {\left\{ {I_{n}}  \right.} \setminus {\left.
{\{i\}} \right\}} $.
\end{lemma}
\begin{lemma}\label{lemma:L_ij} \cite{ky1}
Let $L_{i{\kern 1pt} j} $ be the left $ij$-th cofactor of
 ${\rm {\bf A}}\in {\rm M}\left( {n,{\mathbb{H}}} \right)$, that
 is,
$ {\rm{cdet}} _{{j}}\, {\rm {\bf A}} = {{\sum\limits_{i = 1}^{n}
{L_{i{\kern 1pt} j} \cdot a_{i{\kern 1pt} j}} }}$ for all $j
=1,\ldots,n$. Then
\[
L_{i{\kern 1pt} j} = {\left\{ {\begin{array}{*{20}c}
 -{\rm{cdet}} _{i}\, {\rm {\bf A}}_{i{\kern 1pt} .}^{j{\kern 1pt}j} \left( {{\rm {\bf a}}_{j{\kern 1pt}. } }\right),& {i \ne
j},\\
 {\rm{cdet}} _{k}\, {\rm {\bf A}}^{j\, j},& {i = j},
\\
\end{array} }\right.}
\]
\noindent where  ${\rm {\bf A}}_{i{\kern 1pt} .}^{jj} \left( {{\rm
{\bf a}}_{j{\kern 1pt} .} } \right)$ is obtained from ${\rm {\bf
A}}$
 by replacing the $i$th row with the $j$th row, and then by
deleting both the $j$th row and  column, $k = \min {\left\{
{J_{n}} \right.} \setminus {\left. {\{j\}} \right\}} $.
\end{lemma}

Since these matrix functionals do not satisfy the axioms of
noncommutative determinant, then
 they are called "determinants" shareware.
But by the following theorems we introduce the concepts of a
determinant of a Hermitian matrix and a double determinant  which
both satisfy the axioms of  noncommutative determinant and can  be
expanded  by cofactors along an arbitrary row or column using the
row-column determinants. This enables us to obtain determinantal
representations of the inverse matrix.

\begin{theorem} \cite{ky1}\label{kyrc2}
If ${\rm {\bf A}} = \left( {a_{ij}}  \right) \in {\rm M}\left(
{n,{\rm {\mathbb{H}}}} \right)$ is Hermitian, then ${\rm{rdet}}
_{1} {\rm {\bf A}} = \cdots = {\rm{rdet}} _{n} {\rm {\bf A}} =
{\rm{cdet}} _{1} {\rm {\bf A}} = \cdots = {\rm{cdet}} _{n} {\rm
{\bf A}} \in {\rm {\mathbb{R}}}.$
\end{theorem}
\begin{remark}
\label{remark: determinant of hermitian matrix}
 Since Theorem \ref{kyrc2}, we can
define the determinant of a  Hermitian matrix ${\rm {\bf A}}\in
{\rm M}\left( {n,{\rm {\mathbb{H}}}} \right)$ putting for all $i
=1,\ldots,n$
\[\det {\rm {\bf A}}: = {\rm{rdet}}_{{i}}\,
{\rm {\bf A}} = {\rm{cdet}} _{{i}}\, {\rm {\bf A}}. \]
\end{remark}
\begin{theorem}\cite{ky1} \label{inver_her} If for a Hermitian matrix ${\rm {\bf A}}\in {\rm M}\left( {n,{\rm
{\mathbb{H}}}} \right)$,
\[\det {\rm {\bf A}} \ne 0,\]
then there exist a unique right inverse  matrix $(R{\rm {\bf
A}})^{ - 1}$ and a unique left inverse matrix $(L{\rm {\bf A}})^{
- 1}$ of
 ${\rm {\bf A}}$, where $\left( {R{\rm {\bf A}}} \right)^{ - 1} = \left( {L{\rm {\bf A}}}
\right)^{ - 1} = :{\rm {\bf A}}^{ - 1}$, and they possess the
following determinantal representations

\begin{equation}\label{eq:inver_her_R}
  \left( {R{\rm {\bf A}}} \right)^{ - 1} = {\frac{{1}}{{\det {\rm
{\bf A}}}}}
\begin{pmatrix}
  R_{11} & R_{21} & \cdots & R_{n1}\\
  R_{12} & R_{22} & \cdots & R_{n2}\\
  \cdots & \cdots & \cdots& \cdots\\
  R_{1n} & R_{2n} & \cdots & R_{nn}
\end{pmatrix},
\end{equation}
\begin{equation}\label{eq:inver_her_L}
  \left( {L{\rm {\bf A}}} \right)^{ - 1} = {\frac{{1}}{{\det {\rm
{\bf A}}}}}
\begin{pmatrix}
  L_{11} & L_{21} & \cdots & L_{n1} \\
  L_{12} & L_{22} & \cdots & L_{n2} \\
  \cdots & \cdots & \cdots & \cdots \\
  L_{1n} & L_{2n} & \cdots & L_{nn}
\end{pmatrix},
\end{equation}
 where $R_{ij}$,  $L_{ij}$ are right and left $ij$th cofactors of
${\rm {\bf
 A}}$ respectively for all $ i,j =
1,...,n$.
\end{theorem}
\begin{theorem}\cite{ky1}\label{theorem:detAA}
If ${\rm {\bf A}} \in {\rm M}\left( {n,{\rm {\mathbb{H}}}}
\right)$, then $\det {\rm {\bf A}}{\rm {\bf A}}^{ *} = \det {\rm
{\bf A}}^{ *} {\rm {\bf A}}$.
\end{theorem}
According to Theorem \ref{theorem:detAA} we introduce the concept
of a double determinant.
\begin{definition}
The determinant of the corresponding Hermitian matrices, (${{\rm
{\bf A}}^{ *} {\rm {\bf A}}}$ or  ${{\rm {\bf A}}{\rm {\bf A}}^{
*} }$), of ${\rm {\bf A}} \in {\rm M}\left( {n,{\rm {\mathbb{H}}}}
\right)$
 is called its double determinant, i.e.
${\rm{ddet}}{ \rm{\bf A}}: = \det \left( {{\rm {\bf A}}^{ *} {\rm
{\bf A}}} \right) = \det \left( {{\rm {\bf A}}{\rm {\bf A}}^{ *} }
\right)$.
\end{definition}

Suppose ${\rm {\bf A}}_{}^{i{\kern 1pt} j} $ denotes the submatrix
of ${\rm {\bf A}}$ obtained by deleting both the $i$th row and the
$j$th column. Let ${\rm {\bf a}}_{.j} $ be the $j$th column and
${\rm {\bf a}}_{i.} $ be the $i$th row of ${\rm {\bf A}}$. Denote
by ${\rm {\bf a}}^{\ast}_{.j} $ and ${\rm {\bf a}}^{\ast}_{i.} $
the $j$th column  and the $i$th row of a Hermitian adjoint matrix
${\rm {\bf A}}^{\ast}$ as well. Suppose ${\rm {\bf A}}_{.j} \left(
{{\rm {\bf b}}} \right)$ denotes the matrix obtained from ${\rm
{\bf A}}$ by replacing its $j$th column with the column ${\rm {\bf
b}}$, and ${\rm {\bf A}}_{i.} \left( {{\rm {\bf b}}} \right)$
denotes the matrix obtained from ${\rm {\bf A}}$ by replacing its
$i$th row with the row ${\rm {\bf b}}$.

We have the following theorem on the determinantal representation
of the inverse matrix over ${\mathbb{H}}$.
\begin{theorem}\cite{ky1} \label{theorem:deter_inver} The necessary and sufficient condition of invertibility
of  ${\rm {\bf A}} \in {\rm M}(n,{{\rm {\mathbb{H}}}})$ is
${\rm{ddet}} {\rm {\bf A}} \ne 0$. Then there exists ${\rm {\bf
A}}^{ - 1} = \left( {L{\rm {\bf A}}} \right)^{ - 1} = \left(
{R{\rm {\bf A}}} \right)^{ - 1}$, where
\begin{equation}
\label{eq:det_inv_rdet} \left( {L{\rm {\bf A}}} \right)^{ - 1}
=\left( {{\rm {\bf A}}^{ *}{\rm {\bf A}} } \right)^{ - 1}{\rm {\bf
A}}^{ *} ={\frac{{1}}{{{\rm{ddet}}{ \rm{\bf A}} }}}
\begin{pmatrix}
  {\mathbb{L}} _{11} & {\mathbb{L}} _{21}& \ldots & {\mathbb{L}} _{n1} \\
  {\mathbb{L}} _{12} & {\mathbb{L}} _{22} & \ldots & {\mathbb{L}} _{n2} \\
  \ldots & \ldots & \ldots & \ldots \\
 {\mathbb{L}} _{1n} & {\mathbb{L}} _{2n} & \ldots & {\mathbb{L}} _{nn}
\end{pmatrix},\end{equation}
\begin{equation}
\label{eq:det_inv_cdet}
 \left( {R{\rm {\bf A}}} \right)^{ - 1} = {\rm {\bf
A}}^{ *} \left( {{\rm {\bf A}}{\rm {\bf A}}^{ *} } \right)^{ - 1}
= {\frac{{1}}{{{\rm{ddet}}{ \rm{\bf A}}^{ *} }}}
\begin{pmatrix}
 {\mathbb{R}} _{11} & {\mathbb{R}} _{21} &\ldots & {\mathbb{R}} _{n1} \\
 {\mathbb{R}} _{12} & {\mathbb{R}} _{22} &\ldots & {\mathbb{R}} _{n2}  \\
 \ldots  & \ldots & \ldots & \ldots \\
 {\mathbb{R}} _{1n} & {\mathbb{R}} _{2n} &\ldots & {\mathbb{R}} _{nn}
\end{pmatrix}
\end{equation}
and \[{\mathbb{L}} _{ij} = {\rm{cdet}} _{j} ({\rm {\bf
A}}^{\ast}{\rm {\bf A}})_{.j} \left( {{\rm {\bf a}}_{.{\kern 1pt}
i}^{ *} } \right), \,\,\,{\mathbb{R}} _{\,{\kern 1pt} ij} =
{\rm{rdet}}_{i} ({\rm {\bf A}}{\rm {\bf A}}^{\ast})_{i.} \left(
{{\rm {\bf a}}_{j.}^{ *} }  \right),\] for all $i,j = 1,...,n.$
\end{theorem}

\begin{remark}
 Since by Theorem \ref{theorem:deter_inver} \[{\rm{ddet}} {\rm {\bf A}} =
{\rm{cdet}} _{j} \left( {{\rm {\bf A}}^{ *} {\rm {\bf A}}} \right)
= {\sum\limits_{i} {{\mathbb{L}} _{ij} \cdot a_{ij}} },
\,\,{\rm{ddet}} {\rm {\bf A}} = {\rm{rdet}}_{i} \left( {{\rm {\bf
A}}{\rm {\bf A}}^{ *} } \right) = {\sum\limits_{j} {a_{ij} \cdot}
} {\mathbb{R}} _{ i{\kern 1pt}j}
\]
for all $j =1,\ldots,n $, then ${\mathbb{L}} _{ij} $ is called the
left double $ij$th cofactor and ${\mathbb{R}} _{ i{\kern1pt}j} $
is called the right double $ij$th cofactor for the entry $a_{ij}$
of ${\rm {\bf A}} \in {\rm M}\left( {n,{\rm {\mathbb{H}}}}
\right)$.
\end{remark}
\section{ Determinantal representation of the Moore-Penrose
inverse.}

For a quaternion matrix ${\rm {\bf A}} \in {\rm
{\mathbb{H}}}^{m\times n}$, a generalized inverse of ${\rm {\bf
A}}$ is a quaternion matrix ${\rm {\bf X}}$ with following Penrose
conditions
\begin{equation}\label{eq:MP_prop}
\begin{array}{l}
  1)\,\, \left( {{\rm {\bf A}}{\rm {\bf
A}}^{ +} } \right)^{ *}  = {\rm {\bf A}}{\rm {\bf A}}^{ +};\\
  2)\,\, \left( {{\rm {\bf A}}^{ +} {\rm {\bf A}}} \right)^{ *}  = {\rm
{\bf A}}^{ +} {\rm {\bf A}};\\
  3)\,\, {\rm {\bf A}}{\rm {\bf A}}^{ +}
{\rm {\bf A}} = {\rm {\bf A}};\\
  4)\,\,{\rm {\bf A}}^{ +} {\rm {\bf
A}}{\rm {\bf A}}^{ +}  = {\rm {\bf A}}^{ +}.
\end{array}
\end{equation}
\begin{definition}
For ${\rm {\bf A}} \in {\rm {\mathbb{H}}}^{m\times n}$, ${\rm {\bf
X}}\in {\rm {\mathbb{H}}}^{n\times m}$ is said to be a $(i,j,...)$
generalized inverse of ${\rm {\bf A}}$ if ${\rm {\bf X}}$
satisfies Penrose conditions $(i), (j),...$ in (\ref{eq:MP_prop}).
We denote the ${\rm {\bf X}}$ by ${\rm {\bf A}}^{(i,j,...)}$ and
the set of all ${\rm {\bf A}}^{(i,j,...)}$ by ${\rm {\bf
A}}\{i,j,...\}$.
\end{definition}

By \cite{hu} we know that for ${\rm {\bf A}} \in {\rm
{\mathbb{H}}}^{m\times n}$, the ${\rm {\bf A}}\{i,j,...\}$ exists
and the ${\rm {\bf A}}^{(1,2,3,4)}$ exists uniquely. The matrix
${\rm {\bf A}}^{(1,2,3,4)}$ is  called the Moore–Penrose
 inverse of ${\rm {\bf A}}$ and denote ${\rm {\bf A}}^{ +} :={\rm {\bf A}}^{(1,2,3,4)}$.

In \cite{ky3} the determinantal representations of the
Moore-Penrose inverse over the quaternion skew field was derived
based on the limit representation.
\begin{theorem}
If ${\rm {\bf A}} \in {\rm {\mathbb{H}}}^{m\times n}$ and ${\rm
{\bf A}}^{ +}$ is its Moore-Penrose  inverse, then ${\rm {\bf
A}}^{ +}  = {\mathop {\lim} \limits_{\alpha \to 0}} {\rm {\bf
A}}^{ * }\left( {{\rm {\bf A}}{\rm {\bf A}}^{ *}  + \alpha {\rm
{\bf I}}} \right)^{ - 1} = {\mathop {\lim} \limits_{\alpha \to 0}}
\left( {{\rm {\bf A}}^{
* }{\rm {\bf A}} + \alpha {\rm {\bf I}}} \right)^{ - 1}{\rm {\bf
A}}^{ *} $, where $\alpha \in {\rm {\mathbb {R}}}_{ +}  $.
\end{theorem}
\newtheorem{corollary}{Corollary}[section]
\begin{corollary}\label{cor:A+A-1}  If
${\rm {\bf A}}\in {\mathbb H}^{m\times n} $, then the following
statements are true.
 \begin{itemize}
\item [ i)] If $\rm{rank}\,{\rm {\bf A}} = n$, then ${\rm {\bf A}}^{ +}
= \left( {{\rm {\bf A}}^{ *} {\rm {\bf A}}} \right)^{ - 1}{\rm
{\bf A}}^{ * }$ .
\item [ ii)] If $\rm{rank}\,{\rm {\bf A}} =
m$, then ${\rm {\bf A}}^{ +}  = {\rm {\bf A}}^{ * }\left( {{\rm
{\bf A}}{\rm {\bf A}}^{ *} } \right)^{ - 1}.$
\item [ iii)] If $\rm{rank}\,{\rm {\bf A}} = n = m$, then ${\rm {\bf
A}}^{ +}  = {\rm {\bf A}}^{ - 1}$ .
\end{itemize}
\end{corollary}
We shall use the following notations. Let $\alpha : = \left\{
{\alpha _{1} ,\ldots ,\alpha _{k}} \right\} \subseteq {\left\{
{1,\ldots ,m} \right\}}$ and $\beta : = \left\{ {\beta _{1}
,\ldots ,\beta _{k}} \right\} \subseteq {\left\{ {1,\ldots ,n}
\right\}}$ be subsets of the order $1 \le k \le \min {\left\{
{m,n} \right\}}$. By ${\rm {\bf A}}_{\beta} ^{\alpha} $ denote the
submatrix of ${\rm {\bf A}}$ determined by the rows indexed by
$\alpha$ and the columns indexed by $\beta$. Then ${\rm {\bf
A}}{\kern 1pt}_{\alpha} ^{\alpha}$ denotes the principal submatrix
determined by the rows and columns indexed by $\alpha$.
 If ${\rm {\bf A}} \in {\rm
M}\left( {n,{\rm {\mathbb{H}}}} \right)$ is Hermitian, then by
${\left| {{\rm {\bf A}}_{\alpha} ^{\alpha} } \right|}$ denote the
corresponding principal minor of $\det {\rm {\bf A}}$.
 For $1 \leq k\leq n$, denote by $\textsl{L}_{ k,
n}: = {\left\{ {\,\alpha :\alpha = \left( {\alpha _{1} ,\ldots
,\alpha _{k}} \right),\,{\kern 1pt} 1 \le \alpha _{1} \le \ldots
\le \alpha _{k} \le n} \right\}}$ the collection of strictly
increasing sequences of $k$ integers chosen from $\left\{
{1,\ldots ,n} \right\}$. For fixed $i \in \alpha $ and $j \in
\beta $, let $I_{r,\,m} {\left\{ {i} \right\}}: = {\left\{
{\,\alpha :\alpha \in L_{r,m} ,i \in \alpha}  \right\}}{\rm ,}
\quad J_{r,\,n} {\left\{ {j} \right\}}: = {\left\{ {\,\beta :\beta
\in L_{r,n} ,j \in \beta}  \right\}}$.

The following theorem and remarks introduce the determinantal
representations of the Moore-Penrose inverse which we shall use
below.
\begin{theorem}\label{theor:det_repr_MP} \cite{ky3}
If ${\rm {\bf A}} \in {\rm {\mathbb{H}}}_{r}^{m\times n} $, then
the Moore-Penrose inverse  ${\rm {\bf A}}^{ +} = \left( {a_{ij}^{
+} } \right) \in {\rm {\mathbb{H}}}_{}^{n\times m} $ possess the
following determinantal representations:
\begin{equation}
\label{eq:det_repr_A*A}
 a_{ij}^{ +}  = {\frac{{{\sum\limits_{\beta
\in J_{r,\,n} {\left\{ {i} \right\}}} {{\rm{cdet}} _{i} \left(
{\left( {{\rm {\bf A}}^{ *} {\rm {\bf A}}} \right)_{\,. \,i}
\left( {{\rm {\bf a}}_{.j}^{ *} }  \right)} \right){\kern 1pt}
{\kern 1pt} _{\beta} ^{\beta} } } }}{{{\sum\limits_{\beta \in
J_{r,\,\,n}} {{\left| {\left( {{\rm {\bf A}}^{ *} {\rm {\bf A}}}
\right){\kern 1pt} _{\beta} ^{\beta} }  \right|}}} }}},
\end{equation}
or
\begin{equation}
\label{eq:det_repr_AA*} a_{ij}^{ +}  =
{\frac{{{\sum\limits_{\alpha \in I_{r,m} {\left\{ {j} \right\}}}
{{\rm{rdet}} _{j} \left( {({\rm {\bf A}}{\rm {\bf A}}^{ *}
)_{j\,.\,} ({\rm {\bf a}}_{i.\,}^{ *} )} \right)\,_{\alpha}
^{\alpha} } }}}{{{\sum\limits_{\alpha \in I_{r,\,m}}  {{\left|
{\left( {{\rm {\bf A}}{\rm {\bf A}}^{ *} } \right){\kern 1pt}
_{\alpha} ^{\alpha} } \right|}}} }}}.
\end{equation}
for all $i = 1,...,n $, $j =1,...,m $.
\end{theorem}
\begin{remark}\label{rem:det_repr_MP_A*A}
If ${\rm rank}\,{\rm {\bf A}} = n$, then by Corollary
\ref{cor:A+A-1} ${\rm {\bf A}}^{ +}  = \left( {{\rm {\bf A}}^{ *}
{\rm {\bf A}}} \right)^{ - 1}{\rm {\bf A}}^{ *} $. Considering
$\left( {{\rm {\bf A}}^{ *} {\rm {\bf A}}} \right)^{ - 1}$ as a
left inverse, we get the following representation of  ${\rm {\bf
A}}^{ +} $:
\begin{equation}
\label{eq:det_repr_MP_A*A}{\rm {\bf A}}^{ +}  = {\frac{{1}}{{\rm
{ddet} {\rm {\bf A}}}}}
\begin{pmatrix}
  {{\rm{cdet}} _{1} ({{\rm {\bf A}}^{ *} {\rm {\bf
A}}})_{.\,1} \left( {{\rm {\bf a}}_{.\,1}^{ *} }  \right)} &
\ldots & {{\rm{cdet}} _{1} ({{\rm {\bf A}}^{ *} {\rm {\bf
A}}})_{.\,1} \left(
{{\rm {\bf a}}_{.\,m}^{ *} }  \right)} \\
  \ldots & \ldots & \ldots \\
  {{\rm{cdet}} _{n} ({{\rm {\bf A}}^{ *} {\rm {\bf A}}})_{.\,n}
 \left( {{\rm {\bf a}}_{.\,1}^{ *} }  \right)} & \ldots &
 {{\rm{cdet}} _{n} ({{\rm {\bf A}}^{ *} {\rm {\bf
A}}})_{.\,n} \left( {{\rm {\bf a}}_{.\,m}^{ *} }  \right)}.
\end{pmatrix}
\end{equation}
If $m > n$, then by Theorem \ref{theor:det_repr_MP} for ${\rm {\bf
A}}^{ +} $ we have (\ref{eq:det_repr_A*A}) as well.
\end{remark}
\begin{remark}\label{rem:det_repr_MP_AA*}
If ${\rm rank}\,{\rm {\bf A}} = m$, then by Corollary
\ref{cor:A+A-1} ${\rm {\bf A}}^{ +}  = {\rm {\bf A}}^{ *} \left(
{{\rm {\bf A}}{\rm {\bf A}}^{ *} } \right)^{ - 1}$. Considering
$\left( {{\rm {\bf A}}{\rm {\bf A}}^{ *} } \right)^{ - 1}$ as a
right inverse, we get the following representation of ${\rm {\bf
A}}^{ +} $:
\begin{equation}
\label{kyr9} {\rm {\bf A}}^{ +}  = {\frac{{1}}{{\rm {ddet} {{\rm
{\bf A}} }}}}
\begin{pmatrix}
 {{\rm{rdet}} _{1} ({\rm {\bf A}} {\rm {\bf A}}^{*})_{1.} \left( {{\rm {\bf a}}_{1.}^{ *} }  \right)} & \ldots &
 {{\rm{rdet}} _{m} ({\rm {\bf A}} {\rm {\bf A}}^{*})_{m.} \left( {{\rm {\bf a}}_{1.}^{ *} }  \right)}\\
 \ldots & \ldots & \ldots \\
  {{\rm{rdet}} _{1} ({\rm {\bf A}} {\rm {\bf A}}^{*})_{1.} \left( {{\rm {\bf a}}_{n.}^{ *} }  \right)} & \ldots &
  {{\rm{rdet}} _{m} ({\rm {\bf A}} {\rm {\bf A}}^{*})_{m\,.}
\left( {{\rm {\bf a}}_{n\,.}^{ *} }  \right)}
\end{pmatrix}.
\end{equation}
If $m < n$, then by Theorem \ref{theor:det_repr_MP} for ${\rm {\bf
A}}^{ +} $ we also have (\ref{eq:det_repr_AA*}).
\end{remark}
\section{ Explicit representation formulas for the minimum norm least squares solutions of some  quaternion
matrix equations }

Denote by $\|{\rm {\bf A}}\|$ the Frobenius norm of the quaternion
matrix ${\rm {\bf A}}\in {\mathbb{H}}^{m\times n}$.
\begin{definition}
Consider a  matrix equation
\begin{equation}\label{eq:AX=B}
 {\rm {\bf A}}{\rm {\bf X}} = {\rm {\bf B}},
\end{equation}
\noindent where ${\rm {\bf A}}\in {\mathbb{H}}^{m\times n},{\rm
{\bf B}} \in {\mathbb{H}}^{m\times s} $ are given, ${\rm {\bf X}}
\in {\mathbb{H}}^{n\times s}$ is unknown. Suppose
\[H_{R}=\{{\rm {\bf X}}|{\rm {\bf X}} \in
{\mathbb{H}}^{n\times s}, \|{\rm {\bf A}}{\rm {\bf X}} - {\rm {\bf
B}}\|=\min \},\]  Then matrices ${\rm {\bf X}} \in
{\mathbb{H}}^{n\times s}$ such that ${\rm {\bf X}}\in H_{R}$ are
called least squares solutions of the matrix equation
(\ref{eq:AX=B}).
 If
${\rm {\bf X}}_{LS}={\min}_ {{\rm {\bf X}}\in H_{R}}\|{\rm {\bf
X}}\|$, then ${\rm {\bf X}}_{LS}$ is called the minimum norm least
squares solution of (\ref{eq:AX=B}).
\end{definition}
If  the equation (\ref{eq:AX=B}) has no precision solutions, then
${\rm {\bf X}}_{LS}$ is its optimal approximation.

The following important theorem is well-known.
\begin{theorem}(\cite{ji2})\label{theor:LS} The least squares solutions of (\ref{eq:AX=B}) are
\[{\rm {\bf X}} = {\rm {\bf A}}^{+}{\rm {\bf B}} + ({\rm {\bf
I}}_{n} - {\rm {\bf A}}^{+} {\rm {\bf A}}){\rm {\bf C}},\] where
${\rm {\bf C}}\in {\mathbb{H}}^{n\times s}$ is an arbitrary
quaternion matrix and the minimum norm least squares solution is
${\rm {\bf X}}_{LS}={\rm {\bf A}}^{+}{\rm {\bf B}}$.
\end{theorem}
 We denote ${\rm {\bf A}}^{ \ast}{\rm {\bf B}}=:\hat{{\rm
{\bf B}}}= (\hat{b}_{ij})\in {\mathbb{H}}^{n\times s}$.
\begin{theorem}
\begin{enumerate}
\item[(i)]  If ${\rm rank}\,{\rm {\bf A}} = r \le m < n$, then for the minimum norm
least square least squares solution ${\rm {\bf
X}}_{LS}=(x_{ij})\in {\mathbb{H}}^{n\times s}$ of (\ref{eq:AX=B})
for all $i =1,...,n $, $j =1,...,s $, we have
\begin{equation}
\label{eq:ls_AX} x_{ij} = {\frac{{{\sum\limits_{\beta \in
J_{r,\,n} {\left\{ {i} \right\}}} {{\rm{cdet}} _{i} \left( {\left(
{{\rm {\bf A}}^{ *} {\rm {\bf A}}} \right)_{\,.\,i} \left(
{\hat{{\rm {\bf b}}}_{.j}} \right)} \right){\kern 1pt} {\kern 1pt}
_{\beta} ^{\beta} } } }}{{{\sum\limits_{\beta \in J_{r,\,\,n}}
{{\left| {\left( {{\rm {\bf A}}^{ *} {\rm {\bf A}}} \right){\kern
1pt} {\kern 1pt} _{\beta} ^{\beta} }  \right|}}} }}}.
\end{equation}
\item[(ii)] If ${\rm rank}\,{\rm {\bf A}} = n$, then  for  ${\rm {\bf
X}}_{LS}=(x_{ij})\in {\mathbb{H}}^{n\times s}$ of (\ref{eq:AX=B})
for all $i = 1,...,n $, $j =1,...,s $, we have
\begin{equation}
\label{eq:ls_AX_full} x_{i\,j} = {\frac{{{\rm cdet} _{i} ({\rm
{\bf A}}^{ *} {\rm {\bf A}})_{.\,i} \left( {\hat{{\rm {\bf
b}}}_{.j}} \right)}}{{ {\rm ddet} {\rm {\bf A}}}}}
\end{equation}
\noindent where $\hat{{\rm {\bf b}}}_{.j}$ is the $j$th column of
$\hat{{\rm {\bf B}}}$ for all $j = 1,...,s $.
\end{enumerate}
\end{theorem}
{\textit{Proof.}}
 (i) If ${\rm rank}\,{\rm {\bf A}} = r
\le m < n$, then by Theorem  \ref{theor:det_repr_MP} we can
represent the matrix ${\rm {\bf A}}^{ +} $ by
(\ref{eq:det_repr_A*A}). Therefore, we obtain for all $i = 1,...,n
$, $j =1,...,s $
\[
x_{ij}  =\sum_{k=1}^{m}
a_{ik}^{+}b_{kj}=\sum_{k=1}^{m}{\frac{{{\sum\limits_{\beta \in
J_{r,\,n} {\left\{ {i} \right\}}} {{\rm{cdet}} _{i} \left( {\left(
{{\rm {\bf A}}^{ *} {\rm {\bf A}}} \right)_{\,. \,i} \left( {{\rm
{\bf a}}_{.k}^{ *} }  \right)} \right) {\kern 1pt} _{\beta}
^{\beta} } } }}{{{\sum\limits_{\beta \in J_{r,\,\,n}} {{\left|
{\left( {{\rm {\bf A}}^{ *} {\rm {\bf A}}} \right){\kern 1pt}
_{\beta} ^{\beta} }  \right|}}} }}}\cdot b_{kj}=
\]
\[
{\frac{{{\sum\limits_{\beta \in J_{r,\,n} {\left\{ {i}
\right\}}}\sum_{k=1}^{m} {{\rm{cdet}} _{i} \left( {\left( {{\rm
{\bf A}}^{ *} {\rm {\bf A}}} \right)_{\,. \,i} \left( {{\rm {\bf
a}}_{.k}^{ *} }  \right)} \right) {\kern 1pt} _{\beta} ^{\beta} }
} }\cdot b_{kj}} {{{\sum\limits_{\beta \in J_{r,\,\,n}} {{\left|
{\left( {{\rm {\bf A}}^{ *} {\rm {\bf A}}} \right){\kern 1pt}
_{\beta} ^{\beta} }  \right|}}} }}}
\]
Since ${\sum\limits_{k} {{\rm {\bf a}}_{.\,k}^{ *}  b_{kj}} }=
\left( {{\begin{array}{*{20}c}
 {{\sum\limits_{k} {a_{1k}^{ *}  b_{kj}} } } \hfill \\
 {{\sum\limits_{k} {a_{2k}^{ *}  b_{kj}} } } \hfill \\
 { \vdots}  \hfill \\
 {{\sum\limits_{k} {a_{nk}^{ *}  b_{kj}} } } \hfill \\
\end{array}} } \right) = \hat{{\rm {\bf b}}}_{.j}$ and denoting the $j$th
column of $\hat{{\rm {\bf B}}}$ by $\hat{{\rm {\bf b}}}_{.j} $ as
well, then it follows (\ref{eq:ls_AX}).

\noindent ii) If ${\rm rank}\,{\rm {\bf A}} = n$,  then by
Corollary \ref{cor:A+A-1}, ${\rm {\bf A}}^{ +}  = \left( {{\rm
{\bf A}}^{ *} {\rm {\bf A}}} \right)^{ - 1}{\rm {\bf A}}^{ *} $.
Representing $\left( {{\rm {\bf A}}^{ *} {\rm {\bf A}}} \right)^{
- 1}$ by (\ref{eq:inver_her_L}), we obtain  for all $i =1,...,n $,
$j =1,...,s $
\[
x_{ij}  = {\frac{{1}}{{{\rm ddet} {\rm {\bf A}}}}}{\sum\limits_{k
= 1}^{n} {{L}_{ki} \hat{b}_{kj}} },
\]
\noindent where ${L}_{ij}$ is a left $ij$th
 cofactor of $({\rm {\bf A}}^{ \ast}{\rm {\bf A}})$ for all $i,j = 1,...,n
$. From this by  Lemma \ref{lemma:L_ij}  and denoting the $j$th
column of $\hat{{\rm {\bf B}}}$ by $\hat{{\rm {\bf b}}}_{.j} $, it
follows (\ref{eq:ls_AX_full}). $\blacksquare$
 \begin{corollary}(Theorem 3.1 in \cite{ky2} )
 Suppose
\begin{equation}\label{kyr5}
 {\rm {\bf A}}{\rm {\bf X}} = {\rm {\bf B}}
\end{equation}
\noindent is a right matrix equation, where ${\left\{ {{\rm {\bf
A}},{\rm {\bf B}}} \right\}} \in {\rm M}(n,{\rm {\mathbb{H}}} )$
are given, ${\rm {\bf X}} \in {\rm M}(n,{\rm {\mathbb{H}}} )$ is
unknown. If ${\rm ddet} {\rm {\bf A}} \ne 0$, then (\ref{kyr5})
has a unique solution, and the solution is
\begin{equation}
\label{kyr6} x_{i\,j} = {\frac{{{\rm cdet} _{i} ({\rm {\bf A}}^{
*} {\rm {\bf A}})_{.\,i} \left( {\hat{{\rm {\bf b}}}_{.j}}
\right)}}{{ {\rm ddet} {\rm {\bf A}}}}}
\end{equation}
\noindent where $\hat{{\rm {\bf b}}}_{.j}$ is the $j$th column of
$\hat{{\rm {\bf B}}}$ for all $i,j = 1,...,n $.
 \end{corollary}
\begin{definition}
Consider a  matrix equation
\begin{equation}\label{eq:XA=B}
 {\rm {\bf X}}{\rm {\bf A}} = {\rm {\bf B}},
\end{equation}
\noindent  where ${\rm {\bf A}}\in {\mathbb{H}}^{m\times n},{\rm
{\bf B}} \in {\mathbb{H}}^{s\times n} $ are given, ${\rm {\bf X}}
\in {\mathbb{H}}^{s\times m}$ is unknown. Suppose
\[H_{L}=\{{\rm {\bf X}}|\,{\rm {\bf X}} \in
{\mathbb{H}}^{s\times m}, \|{\rm {\bf X}}{\rm {\bf A}} - {\rm {\bf
B}}\|=\min \}.\] Then matrices ${\rm {\bf X}} \in
{\mathbb{H}}^{s\times m}$ such that ${\rm {\bf X}}\in H_{L}$ are
called least squares solutions of the matrix equation
(\ref{eq:XA=B}).
 If
${\rm {\bf X}}_{LS}={\min}_ {{\rm {\bf X}}\in H_{L}}\|{\rm {\bf
X}}\|$, then ${\rm {\bf X}}_{LS}$ is called the minimum norm least
squares solution of (\ref{eq:XA=B}).
\end{definition}
The following  theorem can be obtained by analogy to Theorem
\ref{theor:LS}.
\begin{theorem} The least squares solutions of (\ref{eq:XA=B}) are
\[{\rm {\bf X}} = {\rm {\bf B}}{\rm {\bf A}}^{+} + ({\rm {\bf
I}}_{m} -  {\rm {\bf A}}{\rm {\bf A}}^{+}){\rm {\bf C}},\] where
${\rm {\bf C}}\in {\mathbb{H}}^{n\times s}$ is an arbitrary
quaternion matrix and the minimum norm least squares solution is
${\rm {\bf X}}_{LS}={\rm {\bf B}}{\rm {\bf A}}^{+}$.
\end{theorem}
 We denote  ${\rm {\bf B}}{\rm {\bf A}}^{
\ast}=:\check{{\rm {\bf B}}}= (\check{b}_{ij})\in
{\mathbb{H}}^{s\times m}$.
\begin{theorem}
\begin{enumerate}
\item[(i)] If ${\rm rank}\,{\rm {\bf A}} = r \le n < m$, then
for the minimum norm least squares solution ${\rm {\bf
X}}_{LS}=(x_{ij})\in {\mathbb{H}}^{s\times m}$ of (\ref{eq:XA=B})
for all $i = 1,...,s $, $j = 1,...,m $, we have
\begin{equation}
\label{eq:ls_XA} x_{ij} = {\frac{{{\sum\limits_{\alpha \in I_{r,m}
{\left\{ {j} \right\}}} {{\rm{rdet}} _{j} \left( {\left( {{\rm
{\bf A}}{\rm {\bf A}}^{ *} } \right)_{\,j\,.} \left( {\check{{\rm
{\bf b}}}_{i\,.}} \right)} \right)\,_{\alpha} ^{\alpha} }
}}}{{{\sum\limits_{\alpha \in I_{r,\,m}}  {{\left| {\left( {{\rm
{\bf A}}{\rm {\bf A}}^{ *} } \right) {\kern 1pt} _{\alpha}
^{\alpha} } \right|}}} }}}.
\end{equation}
\item[(ii)] If ${\rm rank}\,{\rm {\bf A}} = m$, then for  ${\rm {\bf
X}}_{LS}=(x_{ij})\in {\mathbb{H}}^{s\times m}$ of (\ref{eq:XA=B})
for all $i = 1,...,s $, $j =1,...,m $, we have
\begin{equation}
\label{eq:ls_XA_full} x_{i\,j} = {\frac{{{\rm rdet} _{j} ({\rm
{\bf A}}{\rm {\bf A}}^{ *} )_{j.\,} \left( {\check{{\rm {\bf
b}}}_{i\,.}} \right)}}{{{\rm ddet} {\rm {\bf A}}}}}
\end{equation}
\noindent where $\check{{\rm {\bf b}}}_{i.}$ is the $i$th row of
$\check{{\rm {\bf B}}}$ for all $i =1,...,s $.
\end{enumerate}
\end{theorem}
{\textit{Proof.}} (i) If ${\rm rank}\,{\rm {\bf A}} = r \le n <m
$, then by Theorem  \ref{theor:det_repr_MP} we can represent the
matrix ${\rm {\bf A}}^{ +} $ by (\ref{eq:det_repr_AA*}).
Therefore, for all $i =1,...,s $, $j =1,...,m $, we get
\[
x_{ij}  =\sum_{k=1}^{m}b_{ik} a_{kj}^{+}=\sum_{k=1}^{m}b_{ik}\cdot
{\frac{{{\sum\limits_{\alpha \in I_{r,m} {\left\{ {j} \right\}}}
{{\rm{rdet}} _{j} \left( {\left( {{\rm {\bf A}}{\rm {\bf A}}^{ *}
} \right)_{\,j\,.} \left( {\bf a}_{k\,.}^{*} \right)}
\right)\,_{\alpha} ^{\alpha} } }}}{{{\sum\limits_{\alpha \in
I_{r,\,m}}  {{\left| {\left( {{\rm {\bf A}}{\rm {\bf A}}^{ *} }
\right) {\kern 1pt} _{\alpha} ^{\alpha} } \right|}}} }}} =
\]
\[
{\frac{{{\sum_{k=1}^{m}b_{ik} \sum\limits_{\alpha \in I_{r,m}
{\left\{ {j} \right\}}} {{\rm{rdet}} _{j} \left( {\left( {{\rm
{\bf A}}{\rm {\bf A}}^{ *} } \right)_{\,j\,.} \left( {\bf
a}_{k\,.}^{*} \right)} \right)\,_{\alpha} ^{\alpha} } } }}
{{{\sum\limits_{\alpha \in I_{r,\,m}}  {{\left| {\left( {{\rm {\bf
A}}{\rm {\bf A}}^{ *} } \right) {\kern 1pt} _{\alpha} ^{\alpha} }
\right|}}} }}}
\]
Since ${\sum\limits_{k} {{  b_{ik}\rm {\bf a}}_{k\,.}^{ *}}
}=\begin{pmatrix}
  \sum\limits_{k} {b_{ik}a_{k1}^{ *} } & \sum\limits_{k} {b_{ik}a_{k2}^{ *} } & \cdots & \sum\limits_{k} {b_{ik}a_{kn}^{ *} }
\end{pmatrix}
 = \check{{\rm {\bf b}}}_{i.}$ and denoting the $i$th row of
$\check{{\rm {\bf B}}}$ by $\check{{\rm {\bf b}}}_{i.}$ as well,
then it follows (\ref{eq:ls_XA}).

\noindent (ii) If ${\rm rank}\,{\rm {\bf A}} = m$,  then by
Corollary \ref{cor:A+A-1} ${\rm {\bf A}}^{ +}  = {\rm {\bf A}}^{ *
}\left( {{\rm {\bf A}}{\rm {\bf A}}^{ *} } \right)^{ - 1}.$
Representing $\left({\rm {\bf A}} {\rm {\bf A}}^{ *}  \right)^{ -
1}$ by (\ref{eq:inver_her_R}), we obtain  for all $i = 1,...,s $,
$j =1,...,m $,
\[
x_{ij}  = {\frac{{1}}{{{\rm ddet} {\rm {\bf A}}}}}{\sum\limits_{k
= 1}^{n}\check{b}_{ik} {{R}_{jk} } },
\]
\noindent where $R_{ij}$ is a left $ij$th
 cofactor of $({\rm {\bf A}}{\rm {\bf A}}^{ \ast})$ for all $i,j =
1,...,m$. From this by  Lemma \ref{lemma:R_ij}  and denoting the
$i$th row of $\check{{\rm {\bf B}}}$ by $\check{{\rm {\bf
b}}}_{i.}$, it follows (\ref{eq:ls_XA_full}). $\blacksquare$
\begin{corollary}(Theorem 3.2 in \cite{ky2} )
 Suppose
\begin{equation}\label{kyr7}
 {\rm {\bf X}}{\rm {\bf A}} = {\rm {\bf B}}
\end{equation}
\noindent is a left matrix equation, where ${\left\{ {{\rm {\bf
A}},{\rm {\bf B}}} \right\}} \in {\rm M}(n,{\rm {\mathbb{H}}} )$
are given, ${\rm {\bf X}} \in {\rm M}(n,{\rm {\mathbb{H}}} )$ is
unknown. If ${\rm ddet} {\rm {\bf A}} \ne 0$, then (\ref{kyr7})
has a unique solution, and the solution is
\begin{equation}
\label{kyr8} x_{i\,j} = {\frac{{{\rm rdet} _{j} ({\rm {\bf A}}{\rm
{\bf A}}^{ *} )_{j.\,} \left( {\check{{\rm {\bf b}}}_{i\,.}}
\right)}}{{{\rm ddet} {\rm {\bf A}}}}}
\end{equation}
\noindent where $\check{{\rm {\bf b}}}_{i.}$ is the $i$th column
of $\check{{\rm {\bf B}}}$ for all $i,j = 1,...,n $.
\end{corollary}
\begin{definition}
Consider a  matrix equation
\begin{equation}\label{eq:AXB=D}
 {\rm {\bf A}}{\rm {\bf X}}{\rm {\bf B}} = {\rm {\bf
D}},
\end{equation}
\noindent  where $ {\rm {\bf A}}\in{\rm {\mathbb{H}}}^{m \times
n},{\rm {\bf B}}\in{\rm {\mathbb{H}}}^{p \times q}, {\rm {\bf
D}}\in{\rm {\mathbb{H}}}^{m \times q}$ are given, $ {\rm {\bf
X}}\in{\rm {\mathbb{H}}}^{n \times p}$ is unknown. Suppose
\[H_{D}=\{{\rm {\bf X}}|\,{\rm {\bf X}} \in
{\mathbb{H}}^{s\times m}, \|{\rm {\bf A}}{\rm {\bf X}}{\rm {\bf
B}} - {\rm {\bf D}}\|=\min \}.\] Then matrices ${\rm {\bf X}} \in
{\mathbb{H}}^{s\times m}$ such that ${\rm {\bf X}}\in H_{D}$ are
called least squares solutions of the matrix equation
(\ref{eq:AXB=D}).
 If
${\rm {\bf X}}_{LS}={\min}_ {{\rm {\bf X}}\in H_{D}}\|{\rm {\bf
X}}\|$, then ${\rm {\bf X}}_{LS}$ is called the minimum norm least
squares solution of (\ref{eq:AXB=D}).
\end{definition}
The following important theorem is well-known.
\begin{theorem}(\cite{we})\label{theor:LS_AXB} The least squares solutions of (\ref{eq:AXB=D}) are
\[{\rm {\bf X}} = {\rm {\bf A}}^{+}{\rm {\bf
D}}{\rm {\bf B}}^{+} + ({\rm {\bf I}}_{n} - {\rm {\bf A}}^{+} {\rm
{\bf A}}){\rm {\bf V}}+{\rm {\bf W}}({\rm {\bf I}}_{p} -{\rm{\bf
B}}{\rm{\bf B}}^{+}),\] where $\{{\rm {\bf V}},{\rm {\bf
W}}\}\subset {\mathbb{H}}^{n\times p}$  are  arbitrary quaternion
matrices and the least squares solution with minimum norm is ${\rm
{\bf X}}_{LS}={\rm {\bf A}}^{+}{\rm {\bf D}}{\rm {\bf B}}^{+}$.
\end{theorem}
 We denote  ${\rm {\bf \widetilde{D}}}= {\rm {\bf
A}}^\ast{\rm {\bf D}}{\rm {\bf B}}^\ast$.
\begin{theorem}\label{theor:AXB=D}
\begin{enumerate}
\item[(i)] If ${\rm rank}\,{\rm {\bf A}} = r_{1} <  m$ and
 ${\rm rank}\,{\rm {\bf B}} = r_{2} < p$, then for the minimum
 norm
least square solution ${\rm {\bf X}}_{LS}=(x_{ij})\in
{\mathbb{H}}^{n\times p}$  of (\ref{eq:AXB=D}) we have
\begin{equation}\label{eq:d^B}
x_{ij} = {\frac{{{\sum\limits_{\beta \in J_{r_{1},\,n} {\left\{
{i} \right\}}} {{\rm{cdet}} _{i} \left( {\left( {{\rm {\bf A}}^{
*} {\rm {\bf A}}} \right)_{\,.\,i} \left( {{{\rm {\bf
d}}}\,_{.\,j}^{{\rm {\bf B}}}} \right)} \right) _{\beta} ^{\beta}
} } }}{{{\sum\limits_{\beta \in J_{r_{1},n}} {{\left| {\left(
{{\rm {\bf A}}^{ *} {\rm {\bf A}}} \right)_{\beta} ^{\beta} }
\right|}} \sum\limits_{\alpha \in I_{r_{2},p}}{{\left| {\left(
{{\rm {\bf B}}{\rm {\bf B}}^{ *} } \right) _{\alpha} ^{\alpha} }
\right|}}} }}},
\end{equation}
or
\begin{equation}\label{eq:d^A}
 x_{ij}={\frac{{{\sum\limits_{\alpha
\in I_{r_{2},q} {\left\{ {j} \right\}}} {{\rm{rdet}} _{j} \left(
{\left( {{\rm {\bf B}}{\rm {\bf B}}^{ *} } \right)_{\,j\,.} \left(
{{{\rm {\bf d}}}\,_{i\,.}^{{\rm {\bf A}}}} \right)}
\right)\,_{\alpha} ^{\alpha} } }}}{{{\sum\limits_{\beta \in
J_{r_{1},p}} {{\left| {\left( {{\rm {\bf A}}^{ *} {\rm {\bf A}}}
\right) _{\beta} ^{\beta} } \right|}}\sum\limits_{\alpha \in
I_{r_{2},q}} {{\left| {\left( {{\rm {\bf B}}{\rm {\bf B}}^{ *} }
\right) _{\alpha} ^{\alpha} } \right|}}} }}},
\end{equation}
where
\begin{equation} \label{eq:def_d^B_m}
   {{{\rm {\bf d}}}_{.\,j}^{{\rm {\bf B}}}}=\left(
\sum\limits_{\alpha \in I_{r_{2},p} {\left\{ {j} \right\}}}
{{\rm{rdet}} _{j} \left( {\left( {{\rm {\bf B}}{\rm {\bf B}}^{ *}
} \right)_{j.} \left( {\tilde{{\rm {\bf d}}}_{k.}} \right)}
\right)_{\alpha} ^{\alpha}} \right)\in{\rm {\mathbb{H}}}^{n \times
1},\,\,\,\,k=1,...,n \end{equation}
\begin{equation} \label{eq:d_A<n}  {{{\rm {\bf d}}}_{i\,.}^{{\rm {\bf A}}}}=\left(
\sum\limits_{\beta \in J_{r_{1},n} {\left\{ {i} \right\}}}
{{\rm{cdet}} _{i} \left( {\left( {{\rm {\bf A}}^{ *}{\rm {\bf A}}
} \right)_{.i} \left( {\tilde{{\rm {\bf d}}}_{.l}} \right)}
\right)_{\beta} ^{\beta}} \right)\in{\rm {\mathbb{H}}}^{1 \times
p},\,\,\,\,l=1,...,p
\end{equation}
 are the column vector and the row vector, respectively.  ${\tilde{{\rm {\bf
 d}}}_{i.}}$ and
${\tilde{{\rm {\bf d}}}_{.j}}$ are the ith row   and the jth
column  of ${\rm {\bf \widetilde{D}}}$ for all $i =1,...,n $, $j
=1,...,p $.

\item[(ii)] If ${\rm rank}\,{\rm {\bf A}} = n$ and ${\rm rank}\,{\rm {\bf B}} = p$, then for  ${\rm {\bf X}}_{LS}=(x_{ij})\in
{\mathbb{H}}^{n\times p}$  of (\ref{eq:AXB=D}) we have
\begin{equation}
\label{eq:AXB_cdetA*A} x_{i\,j} = {\frac{{{\rm cdet} _{i} ({\rm
{\bf A}}^{ *} {\rm {\bf A}})_{.\,i\,} \left( {{\rm {\bf
d}}_{.j}^{{\rm {\bf B}}}} \right)}}{{{\rm ddet} {\rm {\bf A}}\cdot
{\rm ddet} {\rm {\bf B}}}}},
\end{equation}
or
\begin{equation}
\label{eq:AXB_rdetBB*} x_{i\,j} = {\frac{{{\rm rdet} _{j} ({\rm
{\bf B}}{\rm {\bf B}}^{ *} )_{j.\,} \left( {{\rm {\bf
d}}_{i\,.}^{{\rm {\bf A}}}} \right)}}{{{\rm ddet}  {\rm {\bf
A}}\cdot {\rm ddet} {\rm {\bf B}}}}},
\end{equation}
 \noindent where   \begin{equation}\label{eq:d_B_p}{\rm {\bf d}}_{.j}^{{\rm {\bf B}}} : =
\left( {{\rm rdet} _{j} ({\rm {\bf B}}{\rm {\bf B}}^{ *} )_{j.\,}
\left( {\tilde{{\rm {\bf d}}}_{1\,.}^{}} \right),\ldots ,{\rm
rdet} _{j} ({\rm {\bf B}}{\rm {\bf B}}^{ *} )_{j.\,} \left(
{\tilde{{\rm {\bf d}}}_{n\,.}^{}} \right)}
\right)^{T},\end{equation}
\begin{equation}\label{eq:def_d^A}{\rm {\bf d}}_{i\,.}^{{\rm {\bf A}}} : = \left(
{{\rm cdet} _{i} ({\rm {\bf A}}^{ *} {\rm {\bf A}})_{.\,i\,}
\left( {\tilde{{\rm {\bf d}}}_{.1}} \right),\ldots ,{\rm cdet}
_{i} ({\rm {\bf A}}^{ *} {\rm {\bf A}})_{.\,i\,} \left(
{\tilde{{\rm {\bf d}}}_{.n}} \right)} \right)\end{equation} are
respectively the column-vector and the row-vector. $\tilde{{\rm
{\bf d}}}_{i\,.}$  is the  $i$th row of ${\rm {\bf
\widetilde{D}}}$ for all $i =1,...,n $, and $\tilde{{\rm {\bf
d}}}_{.\,j}$ is the $j$th column of ${\rm {\bf \widetilde{D}}}$
for all $j = 1,...,p $.
\item[(iii)]
If  ${\rm rank}\,{\rm {\bf A}} = n$ and  ${\rm rank}\,{\rm {\bf
B}} = r_{2} < p$, then for   ${\rm {\bf X}}_{LS}=(x_{ij})\in
{\mathbb{H}}^{n\times p}$  of (\ref{eq:AXB=D}) we have
\begin{equation}
\label{eq:AXB_detA*A_d^B} x_{ij}={\frac{{{ {{\rm{cdet}} _{i}
\left( {\left( {{\rm {\bf A}}^{ *} {\rm {\bf A}}} \right)_{\,.\,i}
\left( {\widehat{{\rm {\bf d}}}\,_{.\,j}^{{\rm {\bf B}}}} \right)}
\right) } } }}{{{\rm ddet} {\rm {\bf A}}\sum\limits_{\alpha \in
I_{r_{2},p}}{{\left| {\left( {{\rm {\bf B}}{\rm {\bf B}}^{ *} }
\right) _{\alpha} ^{\alpha} } \right|}}} }},
\end{equation}
or
\begin{equation}\label{AXB_detA*A_d^A}
x_{ij}={\frac{{{\sum\limits_{\alpha \in I_{r_{2},p} {\left\{ {j}
\right\}}} {\rm rdet}_{j}{ \left( {\left( {{\rm {\bf B}}{\rm {\bf
B}}^{ *} } \right)_{\,j\,.} \left( {{{\rm {\bf d}}}\,_{i\,.}^{{\rm
{\bf A}}}} \right)} \right)\,_{\alpha} ^{\alpha} } }}}{{{\rm ddet}
{\rm {\bf A}}\sum\limits_{\alpha \in I_{r_{2},p}}{{\left| {\left(
{{\rm {\bf B}}{\rm {\bf B}}^{ *} } \right) _{\alpha} ^{\alpha} }
\right|}}} }},
\end{equation}
where  $ {{\rm {\bf d}}_{.\,j}^{{\rm {\bf B}}}}$ is
(\ref{eq:def_d^B_m})
   and
  ${\rm {\bf
d}}_{i\,.}^{{\rm {\bf A}}}$ is (\ref{eq:def_d^A}).

\item[(iiii)]
If ${\rm rank}\,{\rm {\bf A}} = r_{1} <  n$ and ${\rm rank}\,{\rm
{\bf B}} =  p$, then for   ${\rm {\bf X}}_{LS}=(x_{ij})\in
{\mathbb{H}}^{n\times p}$ of (\ref{eq:AXB=D}) we have
\begin{equation}
\label{eq:AXB_detBB*_d^A} x_{i\,j} = {\frac{{{\rm rdet} _{j} ({\rm
{\bf B}}{\rm {\bf B}}^{ *} )_{j.\,} \left( {{\rm {\bf
d}}_{i\,.}^{{\rm {\bf A}}}} \right)}}{{{\sum\limits_{\beta \in
J_{r_{1},n}} {{\left| {\left( {{\rm {\bf A}}^{ *} {\rm {\bf A}}}
\right) _{\beta} ^{\beta} } \right|}}\cdot {\rm ddet} {\rm {\bf
B}}}}}},
\end{equation}
or
\begin{equation} \label{eq:AXB_detBB*_d^B} x_{i\,j}=
{\frac{{{\sum\limits_{\beta \in J_{r_{1},\,n} {\left\{ {i}
\right\}}} {{\rm{cdet}} _{i} \left( {\left( {{\rm {\bf A}}^{ *}
{\rm {\bf A}}} \right)_{\,.\,i} \left( {{{\rm {\bf
d}}}\,_{.\,j}^{{\rm {\bf B}}}} \right)} \right) _{\beta} ^{\beta}
} } }}{{{\sum\limits_{\beta \in J_{r_{1},n}} {{\left| {\left(
{{\rm {\bf A}}^{ *} {\rm {\bf A}}} \right)_{\beta} ^{\beta} }
\right|}}\cdot{\rm ddet}{\rm {\bf B}}}}}},
\end{equation}
 \noindent where  $ {{\rm {\bf d}}_{.\,j}^{{\rm {\bf B}}}}$ is
(\ref{eq:d_B_p})
   and
 ${\rm {\bf
d}}_{i\,.}^{{\rm {\bf A}}}$ is (\ref{eq:d_A<n}).

\end{enumerate}
\end{theorem}
{\textit{Proof.}} (i) If ${\rm {\bf A}} \in {\rm
{\mathbb{H}}}_{r_{1}}^{m\times n} $, ${\rm {\bf B}} \in {\rm
{\mathbb{H}}}_{r_{2}}^{p\times q} $ and $ r_{1} <  m$, $ r_{2} <
p$, then by Theorem \ref{theor:det_repr_MP} the Moore-Penrose
inverses ${\rm {\bf A}}^{ +} = \left( {a_{ij}^{ +} } \right) \in
{\rm {\mathbb{H}}}_{}^{n\times m} $ and ${\rm {\bf B}}^{ +} =
\left( {a_{ij}^{ +} } \right) \in {\rm {\mathbb{H}}}^{q\times p} $
posses the following determinantal representations respectively,
\[
 a_{ij}^{ +}  = {\frac{{{\sum\limits_{\beta
\in J_{r_{1},\,n} {\left\{ {i} \right\}}} {{\rm{cdet}} _{i} \left(
{\left( {{\rm {\bf A}}^{ *} {\rm {\bf A}}} \right)_{\,. \,i}
\left( {{\rm {\bf a}}_{.j}^{ *} }  \right)} \right){\kern 1pt}
{\kern 1pt} _{\beta} ^{\beta} } } }}{{{\sum\limits_{\beta \in
J_{r_{1},\,n}} {{\left| {\left( {{\rm {\bf A}}^{ *} {\rm {\bf A}}}
\right){\kern 1pt} _{\beta} ^{\beta} }  \right|}}} }}},
\]
\begin{equation}\label{eq:b+}
 b_{ij}^{ +}  =
{\frac{{{\sum\limits_{\alpha \in I_{r_{2},p} {\left\{ {j}
\right\}}} {{\rm{rdet}} _{j} \left( {({\rm {\bf B}}{\rm {\bf B}}^{
*} )_{j\,.\,} ({\rm {\bf b}}_{i.\,}^{ *} )} \right)\,_{\alpha}
^{\alpha} } }}}{{{\sum\limits_{\alpha \in I_{r_{2},p}}  {{\left|
{\left( {{\rm {\bf B}}{\rm {\bf B}}^{ *} } \right){\kern 1pt}
_{\alpha} ^{\alpha} } \right|}}} }}}.
\end{equation}
By Theorem \ref{theor:LS_AXB},  ${\rm {\bf X}}_{LS}={\rm {\bf
A}}^{+}{\rm {\bf D}}{\rm {\bf B}}^{+}$ and entries of ${\rm {\bf
X}}_{LS}=(x_{ij})$ are

\begin{equation}
\label{eq:sum+} x_{ij} = {{\sum\limits_{s = 1}^{q} {\left(
{{\sum\limits_{k = 1}^{m} {{a}_{ik}^{+} d_{ks}} } } \right)}}
{b}_{sj}^{+}}.
\end{equation}
 for all $i=1,...,n$, $j=1,...,p$.

Denote  by $\hat{{\rm {\bf d}}_{.s}}$ the $s$th column of ${\rm
{\bf A}}^{ \ast}{\rm {\bf D}}=:\hat{{\rm {\bf D}}}=
(\hat{d}_{ij})\in {\mathbb{H}}^{n\times q}$ for all $s=1,...,q$.
It follows from ${\sum\limits_{k} { {\rm {\bf a}}_{.\,k}^{
*}}d_{ks} }=\hat{{\rm {\bf d}}_{.\,s}}$ that
\[
\sum\limits_{k = 1}^{m} {{a}_{ik}^{+} d_{ks}}=\sum\limits_{k =
1}^{m}{\frac{{{\sum\limits_{\beta \in J_{r_{1},\,n} {\left\{ {i}
\right\}}} {{\rm{cdet}} _{i} \left( {\left( {{\rm {\bf A}}^{ *}
{\rm {\bf A}}} \right)_{\,. \,i} \left( {{\rm {\bf a}}_{.k}^{ *} }
\right)} \right){\kern 1pt} {\kern 1pt} _{\beta} ^{\beta} } }
}}{{{\sum\limits_{\beta \in J_{r_{1},\,n}} {{\left| {\left( {{\rm
{\bf A}}^{ *} {\rm {\bf A}}} \right){\kern 1pt} _{\beta} ^{\beta}
}  \right|}}} }}}\cdot d_{ks}=
\]
\begin{equation}\label{eq:sum_cdet}
{\frac{{{\sum\limits_{\beta \in J_{r_{1},\,n} {\left\{ {i}
\right\}}}\sum\limits_{k = 1}^{m} {{\rm{cdet}} _{i} \left( {\left(
{{\rm {\bf A}}^{ *} {\rm {\bf A}}} \right)_{\,. \,i} \left( {{\rm
{\bf a}}_{.k}^{ *} } \right)} \right){\kern 1pt} {\kern 1pt}
_{\beta} ^{\beta} } } }\cdot d_{ks}}{{{\sum\limits_{\beta \in
J_{r_{1},\,n}} {{\left| {\left( {{\rm {\bf A}}^{ *} {\rm {\bf A}}}
\right){\kern 1pt} _{\beta} ^{\beta} }  \right|}}}
}}}={\frac{{{\sum\limits_{\beta \in J_{r_{1},\,n} {\left\{ {i}
\right\}}} {{\rm{cdet}} _{i} \left( {\left( {{\rm {\bf A}}^{ *}
{\rm {\bf A}}} \right)_{\,. \,i} \left( \hat{{\rm {\bf d}}_{.\,s}}
\right)} \right){\kern 1pt} {\kern 1pt} _{\beta} ^{\beta} } }
}}{{{\sum\limits_{\beta \in J_{r_{1},\,n}} {{\left| {\left( {{\rm
{\bf A}}^{ *} {\rm {\bf A}}} \right){\kern 1pt} _{\beta} ^{\beta}
}  \right|}}} }}}
\end{equation}
Suppose ${\rm {\bf e}}_{s.}$ and ${\rm {\bf e}}_{.\,s}$ are
respectively the unit row-vector and the unit column-vector whose
components are $0$, except the $s$th components, which are $1$.
Substituting  (\ref{eq:sum_cdet}) and (\ref{eq:b+}) in
(\ref{eq:sum+}), we obtain
\[
x_{ij} =\sum\limits_{s = 1}^{q}{\frac{{{\sum\limits_{\beta \in
J_{r_{1},\,n} {\left\{ {i} \right\}}} {{\rm{cdet}} _{i} \left(
{\left( {{\rm {\bf A}}^{ *} {\rm {\bf A}}} \right)_{\,. \,i}
\left( \hat{{\rm {\bf d}}_{.\,s}} \right)} \right){\kern 1pt}
{\kern 1pt} _{\beta} ^{\beta} } } }}{{{\sum\limits_{\beta \in
J_{r_{1},\,n}} {{\left| {\left( {{\rm {\bf A}}^{ *} {\rm {\bf A}}}
\right){\kern 1pt} _{\beta} ^{\beta} }  \right|}}}
}}}{\frac{{{\sum\limits_{\alpha \in I_{r_{2},p} {\left\{ {j}
\right\}}} {{\rm{rdet}} _{j} \left( {({\rm {\bf B}}{\rm {\bf B}}^{
*} )_{j\,.\,} ({\rm {\bf b}}_{s.\,}^{ *} )} \right)\,_{\alpha}
^{\alpha} } }}}{{{\sum\limits_{\alpha \in I_{r_{2},p}}  {{\left|
{\left( {{\rm {\bf B}}{\rm {\bf B}}^{ *} } \right){\kern 1pt}
_{\alpha} ^{\alpha} } \right|}}} }}}.
\]
Since \begin{equation}\label{eq:prop}\hat{{\rm {\bf
d}}_{.\,s}}=\sum\limits_{l = 1}^{n}{\rm {\bf e}}_{.\,l}\hat{
d_{ls}},\,  {\rm {\bf b}}_{s.\,}^{ *}=\sum\limits_{t =
1}^{p}b_{st}^{*}{\rm {\bf
e}}_{t.},\,\sum\limits_{s=1}^{q}\hat{d_{ls}}b_{st}^{*}=\widetilde{d}_{lt},\end{equation}
then we have
\[
x_{ij} = \]
\[{\frac{{ \sum\limits_{s = 1}^{q}\sum\limits_{t =
1}^{p} \sum\limits_{l = 1}^{n} {\sum\limits_{\beta \in
J_{r_{1},\,n} {\left\{ {i} \right\}}} {{\rm{cdet}} _{i} \left(
{\left( {{\rm {\bf A}}^{ *} {\rm {\bf A}}} \right)_{\,. \,i}
\left( {\rm {\bf e}}_{.\,l} \right)} \right){\kern 1pt} {\kern
1pt} _{\beta} ^{\beta} } } }\hat{
d_{ls}}b_{st}^{*}{\sum\limits_{\alpha \in I_{r_{2},p} {\left\{ {j}
\right\}}} {{\rm{rdet}} _{j} \left( {({\rm {\bf B}}{\rm {\bf B}}^{
*} )_{j\,.\,} ({\rm {\bf e}}_{t.} )} \right)\,_{\alpha} ^{\alpha}
} } }{{{\sum\limits_{\beta \in J_{r_{1},\,n}} {{\left| {\left(
{{\rm {\bf A}}^{ *} {\rm {\bf A}}} \right){\kern 1pt} _{\beta}
^{\beta} }  \right|}}} }{{\sum\limits_{\alpha \in I_{r_{2},p}}
{{\left| {\left( {{\rm {\bf B}}{\rm {\bf B}}^{ *} } \right){\kern
1pt} _{\alpha} ^{\alpha} } \right|}}} }}    }=
\]
\begin{equation}\label{eq:x_ij}
{\frac{{ \sum\limits_{t = 1}^{p} \sum\limits_{l = 1}^{n}
{\sum\limits_{\beta \in J_{r_{1},\,n} {\left\{ {i} \right\}}}
{{\rm{cdet}} _{i} \left( {\left( {{\rm {\bf A}}^{ *} {\rm {\bf
A}}} \right)_{\,. \,i} \left( {\rm {\bf e}}_{.\,l} \right)}
\right){\kern 1pt} {\kern 1pt} _{\beta} ^{\beta} } }
}\,\,\widetilde{d}_{lt}{\sum\limits_{\alpha \in I_{r_{2},p}
{\left\{ {j} \right\}}} {{\rm{rdet}} _{j} \left( {({\rm {\bf
B}}{\rm {\bf B}}^{ *} )_{j\,.\,} ({\rm {\bf e}}_{t.} )}
\right)\,_{\alpha} ^{\alpha} } } }{{{\sum\limits_{\beta \in
J_{r_{1},\,n}} {{\left| {\left( {{\rm {\bf A}}^{ *} {\rm {\bf A}}}
\right){\kern 1pt} _{\beta} ^{\beta} }  \right|}}}
}{{\sum\limits_{\alpha \in I_{r_{2},p}} {{\left| {\left( {{\rm
{\bf B}}{\rm {\bf B}}^{ *} } \right){\kern 1pt} _{\alpha}
^{\alpha} } \right|}}} }}    }.
\end{equation}
Denote by
\[
 d^{{\rm {\bf A}}}_{it}:= \]
\[
{\sum\limits_{\beta \in J_{r_{1},\,n} {\left\{ {i} \right\}}}
{{\rm{cdet}} _{i} \left( {\left( {{\rm {\bf A}}^{ *} {\rm {\bf
A}}} \right)_{\,. \,i} \left( \widetilde{{\rm {\bf d}}}_{.\,t}
\right)} \right){\kern 1pt} _{\beta} ^{\beta} } }= \sum\limits_{l
= 1}^{n} {\sum\limits_{\beta \in J_{r_{1},\,n} {\left\{ {i}
\right\}}} {{\rm{cdet}} _{i} \left( {\left( {{\rm {\bf A}}^{ *}
{\rm {\bf A}}} \right)_{\,. \,i} \left( {\rm {\bf e}}_{.\,l}
\right)} \right){\kern 1pt} {\kern 1pt} _{\beta} ^{\beta} } }
\widetilde{d}_{lt}
\]
the $t$th component  of a row-vector ${\rm {\bf d}}^{{\rm {\bf
A}}}_{i\,.}= (d^{{\rm {\bf A}}}_{i1},...,d^{{\rm {\bf A}}}_{ip})$
for all $t=1,...,p$. Substituting it in (\ref{eq:x_ij}), we have
\[x_{ij} ={\frac{{ \sum\limits_{t = 1}^{p}
 d^{{\rm {\bf A}}}_{it}
}{\sum\limits_{\alpha \in I_{r_{2},p} {\left\{ {j} \right\}}}
{{\rm{rdet}} _{j} \left( {({\rm {\bf B}}{\rm {\bf B}}^{ *}
)_{j\,.\,} ({\rm {\bf e}}_{t.} )} \right)\,_{\alpha} ^{\alpha} } }
}{{{\sum\limits_{\beta \in J_{r_{1},\,n}} {{\left| {\left( {{\rm
{\bf A}}^{ *} {\rm {\bf A}}} \right){\kern 1pt} _{\beta} ^{\beta}
}  \right|}}} }{{\sum\limits_{\alpha \in I_{r_{2},p}} {{\left|
{\left( {{\rm {\bf B}}{\rm {\bf B}}^{ *} } \right){\kern 1pt}
_{\alpha} ^{\alpha} } \right|}}} }}    }.
\]
Since $\sum\limits_{t = 1}^{p}
 d^{{\rm {\bf A}}}_{it}{\rm {\bf e}}_{t.}={\rm {\bf
d}}^{{\rm {\bf A}}}_{i\,.}$, then it follows (\ref{eq:d^A}).

If we denote by
\begin{equation}\label{eq:d^B_den}
 d^{{\rm {\bf B}}}_{lj}:=
\sum\limits_{t = 1}^{p}\widetilde{d}_{lt}{\sum\limits_{\alpha \in
I_{r_{2},p} {\left\{ {j} \right\}}} {{\rm{rdet}} _{j} \left(
{({\rm {\bf B}}{\rm {\bf B}}^{ *} )_{j\,.\,} ({\rm {\bf e}}_{t.}
)} \right)\,_{\alpha} ^{\alpha} } }={\sum\limits_{\alpha \in
I_{r_{2},p} {\left\{ {j} \right\}}} {{\rm{rdet}} _{j} \left(
{({\rm {\bf B}}{\rm {\bf B}}^{ *} )_{j\,.\,} (\widetilde{{\rm {\bf
d}}}_{l.} )} \right)\,_{\alpha} ^{\alpha} } }
\end{equation}

\noindent the $l$th component  of a column-vector ${\rm {\bf
d}}^{{\rm {\bf B}}}_{.\,j}= (d^{{\rm {\bf B}}}_{1j},...,d^{{\rm
{\bf B}}}_{jn})^{T}$ for all $l=1,...,n$ and substitute it in
(\ref{eq:x_ij}), we obtain
\[x_{ij} ={\frac{{  \sum\limits_{l = 1}^{n}
{\sum\limits_{\beta \in J_{r_{1},\,n} {\left\{ {i} \right\}}}
{{\rm{cdet}} _{i} \left( {\left( {{\rm {\bf A}}^{ *} {\rm {\bf
A}}} \right)_{\,. \,i} \left( {\rm {\bf e}}_{.\,l} \right)}
\right){\kern 1pt}  _{\beta} ^{\beta} } } }\,\,d^{{\rm {\bf
B}}}_{lj} }{{{\sum\limits_{\beta \in J_{r_{1},\,n}} {{\left|
{\left( {{\rm {\bf A}}^{ *} {\rm {\bf A}}} \right){\kern 1pt}
_{\beta} ^{\beta} }  \right|}}} }{{\sum\limits_{\alpha \in
I_{r_{2},p}} {{\left| {\left( {{\rm {\bf B}}{\rm {\bf B}}^{ *} }
\right){\kern 1pt} _{\alpha} ^{\alpha} } \right|}}} }}    }.
\]
Since $\sum\limits_{l = 1}^{n}{\rm {\bf e}}_{.l}
 d^{{\rm {\bf B}}}_{lj}={\rm {\bf
d}}^{{\rm {\bf B}}}_{.\,j}$, then it follows (\ref{eq:d^B}).

(ii) If ${\rm rank}\,{\rm {\bf A}} = n$ and ${\rm rank}\,{\rm {\bf
B}} = p$,
  then by Corollary \ref{cor:A+A-1}, ${\rm {\bf
A}}^{ +}  = \left( {{\rm {\bf A}}^{ *} {\rm {\bf A}}} \right)^{ -
1}{\rm {\bf A}}^{ *} $ and ${\rm {\bf B}}^{ +}  = {\rm {\bf B}}^{
* }\left( {{\rm {\bf B}}{\rm {\bf B}}^{ *} } \right)^{ - 1}$.
 If we represent $({\rm {\bf A}}^{ \ast}{\rm {\bf A}})^{ - 1}$ as the left inverse by ( \ref{eq:inver_her_L}) and $ \left( {{\rm {\bf B}}{\rm {\bf B}}^{ *} }
\right)^{ - 1}$ as the right inverse by  (\ref{eq:inver_her_R}) ,
then we obtain
\[\begin{array}{c}
{\rm {\bf X}}_{LS} = ({\rm {\bf A}}^{ \ast}{\rm {\bf A}})^{ -
1}{\rm {\bf A}}^{ \ast}{\rm {\bf D}}{\rm {\bf B}}^{ *} \left(
{{\rm {\bf
B}}{\rm {\bf B}}^{ *} } \right)^{ - 1}=\\
 = \begin{pmatrix}
   x_{11} & x_{12} & \ldots & x_{1p} \\
   x_{21} & x_{22} & \ldots & x_{2p} \\
   \ldots & \ldots & \ldots & \ldots \\
   x_{n1} & x_{n2} & \ldots & x_{np} \
 \end{pmatrix}
 = {\frac{{1}}{{{\rm ddet} {\rm {\bf A}}}}}\begin{pmatrix}
  {L} _{11}^{{\rm {\bf A}}} & {L} _{21}^{{\rm {\bf A}}}&
   \ldots & {L} _{n1}^{{\rm {\bf A}}} \\
  {L} _{12}^{{\rm {\bf A}}} & {L} _{22}^{{\rm {\bf A}}} &
  \ldots & {L} _{n2}^{{\rm {\bf A}}} \\
  \ldots & \ldots & \ldots & \ldots \\
 {L} _{1n}^{{\rm {\bf A}}} & {L} _{2n}^{{\rm {\bf A}}}
  & \ldots & {L} _{nn}^{{\rm {\bf A}}}
\end{pmatrix}\times\\
  \times\begin{pmatrix}
    \tilde{d}_{11} & \tilde{d}_{12} & \ldots & \tilde{d}_{1m} \\
    \tilde{d}_{21} & \tilde{d}_{22} & \ldots & \tilde{d}_{2m} \\
    \ldots & \ldots & \ldots & \ldots \\
    \tilde{d}_{n1} & \tilde{d}_{n2} & \ldots & \tilde{d}_{nm} \
  \end{pmatrix}

{\frac{{1}}{{{\rm ddet} {\rm {\bf A}}}}}
\begin{pmatrix}
 {R} _{\, 11}^{{\rm {\bf B}}} & {R} _{\, 21}^{{\rm {\bf B}}}
 &\ldots & {R} _{\, p1}^{{\rm {\bf B}}} \\
 {R} _{\, 12}^{{\rm {\bf B}}} & {R} _{\, 22}^{{\rm {\bf B}}} &\ldots &
 {R} _{\, p2}^{{\rm {\bf B}}} \\
 \ldots  & \ldots & \ldots & \ldots \\
 {R} _{\, 1p}^{{\rm {\bf B}}} & {R} _{\, 2p}^{{\rm {\bf B}}} &
 \ldots & {R} _{\, pp}^{{\rm {\bf B}}}
\end{pmatrix},
\end{array}
\]
\noindent where ${L}_{ij}^{{\rm {\bf A}}} $ is a left $ij$th
 cofactor of $({\rm {\bf A}}^{ \ast}{\rm {\bf A}})$ for all $i,j =  1,...,n $ and ${R}_{i{\kern 1pt} j}^{{\rm {\bf
B}}} $ is a right  $ij$th cofactor of $\left( {{\rm {\bf B}}{\rm
{\bf B}}^{ *} } \right)$ for all $i,j =1,...,p $. This implies

\begin{equation}
\label{eq:sum} x_{ij} = {\frac{{{\sum\limits_{s = 1}^{p} {\left(
{{\sum\limits_{k = 1}^{n} {{L}_{ki}^{{\rm {\bf A}}}
\tilde{d}_{\,ks}} } } \right)}} {R}_{js}^{{\rm {\bf B}}}} }{{{\rm
ddet} {\rm {\bf A}}\cdot{\rm ddet} {\rm {\bf B}}}}},
\end{equation}
\noindent for all $i =  1,...,n $, $j =1,...,p $. We obtain the
sum in parentheses by Lemma \ref{lemma:L_ij} and denote it,
\[{\sum\limits_{k = 1}^{n}
{{L}_{ki}^{{\rm {\bf A}}} \tilde{d}_{k\,s}} }  = {\rm cdet} _{i}
({\rm {\bf A}}^{ *} {\rm {\bf A}})_{.\,i\,} \left( {\tilde{{\rm
{\bf d}}}_{.\,s}} \right)=:d_{i\,s}^{{\rm {\bf A}}},\]

\noindent where $\tilde{{\rm {\bf d}}}_{.\,s} $ is the $s$th
column-vector of $\tilde{{\rm {\bf D}}}$ for all $s = 1,...,p $.
Suppose  ${{{\rm {\bf d}}}_{i\,.}^{{\rm {\bf A}}}} : = \left(
d_{i\,1}^{{\rm {\bf A}}},\ldots , d_{i\,p}^{{\rm {\bf A}}}
\right)$ is the row-vector for all $i = 1,...,n $. Reducing the
sum ${{\sum\limits_{s = 1}^{n} d_{i\,s}^{{\rm {\bf A}}}}
{R}_{js}^{{\rm {\bf B}}}} $ by Lemma \ref{lemma:R_ij}, we obtain
an analog of Cramer's rule for the minimum norm least squares
solution of (\ref{eq:AXB=D}) by (\ref{eq:AXB_cdetA*A}).

Interchanging the order of summation in (\ref{eq:sum}), we have
\[
x_{ij} = {\frac{{{\sum\limits_{k = 1}^{n} {{ L}_{ki}^{{\rm {\bf
A}}}} }\left( {{\sum\limits_{s = 1}^{p} {\tilde{d}_{\,ks}} }
{R}_{js}^{{\rm {\bf B}}}} \right)}}{{{\rm ddet} {\rm {\bf
A}}\cdot{\rm ddet}{\rm {\bf B}}}}}.\] Further, we obtain the sum
in parentheses by Lemma \ref{lemma:R_ij} and denote it,
 \[{\sum\limits_{s
= 1}^{p} {\tilde{d}_{k\,s}} } {R}_{j\,s}^{{\rm {\bf B}}} = {\rm
rdet} _{j} ({\rm {\bf B}}{\rm {\bf B}}^{ *} )_{j.\,} \left(
{\tilde{{\rm {\bf d}}}_{k\,.}} \right):=d_{k\,j}^{{\rm {\bf
B}}},\] where $\tilde{{\rm {\bf d}}}_{k\,.} $ is the $k$th
row-vector of $\tilde{{\rm {\bf D}}}$ for all $k = 1,...,n $.
Suppose ${\rm {\bf d}}_{.\,j}^{{\rm {\bf B}}} : =
\left(d_{1\,j}^{{\rm {\bf B}}},\ldots ,d_{n\,j}^{{\rm {\bf B}}}
\right)^{T}$
 is the column-vector for all
$j =1,...,p $. Using Lemma \ref{lemma:L_ij} for reducing the sum
${\sum\limits_{k = 1}^{n} {{ L}_{ki}^{{\rm {\bf A}}}}
}d_{k\,j}^{{\rm {\bf B}}} $, we obtain an analog of Cramer's rule
for (\ref{eq:AXB=D}) by (\ref{eq:AXB_rdetBB*}).

(iii) If ${\rm {\bf A}} \in {\rm {\mathbb{H}}}_{r_{1}}^{m\times n}
$, ${\rm {\bf B}} \in {\rm {\mathbb{H}}}_{r_{2}}^{p\times q} $ and
$ r_{1} =n$, $ r_{2} < p$, then by Remark
\ref{rem:det_repr_MP_A*A} and Theorem \ref{theor:det_repr_MP}  the
Moore-Penrose inverses ${\rm {\bf A}}^{ +} = \left( {a_{ij}^{ +} }
\right) \in {\rm {\mathbb{H}}}_{}^{n\times m} $ and ${\rm {\bf
B}}^{ +} = \left( {b_{ij}^{ +} } \right) \in {\rm
{\mathbb{H}}}^{q\times p} $ possess the following determinantal
representations respectively,
\[
 a_{ij}^{ +}  = {\frac{{\rm cdet}_{i} {\left( {{\rm
{\bf A}}^{ *} {\rm {\bf A}}} \right)_{\,. \,i} \left( {{\rm {\bf
a}}_{.j}^{ *} }  \right)}   } {{\rm ddet} { {{\rm {\bf A}} } }}},
\]
\begin{equation}\label{eq:b+2}
 b_{ij}^{ +}  =
{\frac{{{\sum\limits_{\alpha \in I_{r_{2},p} {\left\{ {j}
\right\}}} {\rm rdet}_{j}{ \left( {({\rm {\bf B}}{\rm {\bf B}}^{
*} )_{j\,.\,} ({\rm {\bf b}}_{i.\,}^{ *} )} \right)\,_{\alpha}
^{\alpha} } }}}{{{\sum\limits_{\alpha \in I_{r_{2},p}}  {{\left|
{\left( {{\rm {\bf B}}{\rm {\bf B}}^{ *} } \right){\kern 1pt}
_{\alpha} ^{\alpha} } \right|}}} }}}.
\end{equation}
Since by Theorem \ref{theor:LS_AXB},  ${\rm {\bf X}}_{LS}={\rm
{\bf A}}^{+}{\rm {\bf D}}{\rm {\bf B}}^{+}$, then an entry of
${\rm {\bf X}}_{LS}=(x_{ij})$ is (\ref{eq:sum+}). Denote  by
$\hat{{\rm {\bf d}}_{.s}}$ the $s$th column of ${\rm {\bf A}}^{
\ast}{\rm {\bf D}}=:\hat{{\rm {\bf D}}}= (\hat{d}_{ij})\in
{\mathbb{H}}^{n\times q}$ for all $s=1,...,q$. It follows from
${\sum\limits_{k} { {\rm {\bf a}}_{.\,k}^{ *}}d_{ks} }=\hat{{\rm
{\bf d}}_{.\,s}}$ that
\begin{equation}\label{eq:sum_det}
\sum\limits_{k = 1}^{m} {{a}_{ik}^{+} d_{ks}}=\sum\limits_{k =
1}^{m}{\frac{{\rm cdet}_{i} {\left( {{\rm {\bf A}}^{ *} {\rm {\bf
A}}} \right)_{\,. \,i} \left( {{\rm {\bf a}}_{.\,k}^{ *} }
\right)} } {{\rm ddet} { {{\rm {\bf A}}} }}}\cdot
d_{ks}={\frac{{\rm cdet}_{i} {\left( {{\rm {\bf A}}^{ *} {\rm {\bf
A}}} \right)_{\,. \,i} \left( {\hat{{\rm {\bf d}}_{.\,s}} }
\right)} } {{\rm ddet} { {{\rm {\bf A}}} }}}
\end{equation}
Substituting  (\ref{eq:sum_det}) and (\ref{eq:b+2}) in
(\ref{eq:sum+}), and using (\ref{eq:prop}) we have
\[
x_{ij} =\sum\limits_{s = 1}^{q}{\frac{{\rm cdet}_{i} {\left( {{\rm
{\bf A}}^{ *} {\rm {\bf A}}} \right)_{\,. \,i} \left( {\hat{{\rm
{\bf d}}_{.\,s}} }  \right)} } {{\rm ddet} { {{\rm {\bf A}}}
}}}{\frac{{{\sum\limits_{\alpha \in I_{r_{2},p} {\left\{ {j}
\right\}}} {\rm rdet}_{j}{ \left( {({\rm {\bf B}}{\rm {\bf B}}^{
*} )_{j\,.\,} ({\rm {\bf b}}_{s.\,}^{ *} )} \right)\,_{\alpha}
^{\alpha} } }}}{{{\sum\limits_{\alpha \in I_{r_{2},p}} {{\left|
{\left( {{\rm {\bf B}}{\rm {\bf B}}^{ *} } \right){\kern 1pt}
_{\alpha} ^{\alpha} } \right|}}} }}}=
\]
\[{\frac{{ \sum\limits_{s = 1}^{q}\sum\limits_{t =
1}^{p} \sum\limits_{l = 1}^{n} {\rm cdet}_{i} {\left( {{\rm {\bf
A}}^{ *} {\rm {\bf A}}} \right)_{\,. \,i} \left( {{\rm {\bf
e}}_{.\,l} } \right)} }\hat{ d_{ls}}b_{st}^{*}{\sum\limits_{\alpha
\in I_{r_{2},p} {\left\{ {j} \right\}}}{\rm rdet}_{j} { \left(
{({\rm {\bf B}}{\rm {\bf B}}^{ *} )_{j\,.\,} ({\rm {\bf e}}_{t.}
)} \right)\,_{\alpha} ^{\alpha} } } }{{{{\rm ddet} { {{\rm {\bf
A}}} }} }{{\sum\limits_{\alpha \in I_{r_{2},p}} {{\left| {\left(
{{\rm {\bf B}}{\rm {\bf B}}^{ *} } \right){\kern 1pt} _{\alpha}
^{\alpha} } \right|}}} }}    }=
\]
\begin{equation}\label{eq:x_ij_2}
{\frac{{ \sum\limits_{t = 1}^{p} \sum\limits_{l = 1}^{n} {\rm
cdet}_{i} {\left( {{\rm {\bf A}}^{ *} {\rm {\bf A}}} \right)_{\,.
\,i} \left( {{\rm {\bf e}}_{.\,l} } \right)}
}\,\,\widetilde{d}_{lt}{\sum\limits_{\alpha \in I_{r_{2},p}
{\left\{ {j} \right\}}} {\rm rdet}_{j}{ \left( {({\rm {\bf B}}{\rm
{\bf B}}^{ *} )_{j\,.\,} ({\rm {\bf e}}_{t.} )} \right)\,_{\alpha}
^{\alpha} } } }{{{{\rm ddet} { {{\rm {\bf A}}} }}
}{{\sum\limits_{\alpha \in I_{r_{2},p}} {{\left| {\left( {{\rm
{\bf B}}{\rm {\bf B}}^{ *} } \right){\kern 1pt} _{\alpha}
^{\alpha} } \right|}}} }}    }.
\end{equation}
If we substitute (\ref{eq:d^B_den}) in (\ref{eq:x_ij_2}), then we
get
\[x_{ij} ={\frac{{  \sum\limits_{l = 1}^{n}
{\rm cdet}_{i} {\left( {{\rm {\bf A}}^{ *} {\rm {\bf A}}}
\right)_{\,. \,i} \left( {{\rm {\bf e}}_{.\,l} } \right)}
}\,\,d^{{\rm {\bf B}}}_{lj} }{{{{\rm ddet} { {{\rm {\bf A}}} }}
}{{\sum\limits_{\alpha \in I_{r_{2},p}} {{\left| {\left( {{\rm
{\bf B}}{\rm {\bf B}}^{ *} } \right){\kern 1pt} _{\alpha}
^{\alpha} } \right|}}} }}    }.
\]
Since again $\sum\limits_{l = 1}^{n}{\rm {\bf e}}_{.l}
 d^{{\rm {\bf B}}}_{lj}={\rm {\bf
d}}^{{\rm {\bf B}}}_{.\,j}$, then it follows
(\ref{eq:AXB_detA*A_d^B}), where $ {{\rm {\bf d}}_{.\,j}^{{\rm
{\bf B}}}}$ is (\ref{eq:def_d^B_m}).

If we denote by
\[
 d^{{\rm {\bf A}}}_{it}:= \]
\[
\sum\limits_{l = 1}^{n} {\rm cdet}_{i} {\left( {{\rm {\bf A}}^{ *}
{\rm {\bf A}}} \right)_{\,. \,i} \left( {\widetilde{{\rm {\bf
d}}}_{.\,t} } \right)} = \sum\limits_{l = 1}^{n} {\rm cdet}_{i}
{\left( {{\rm {\bf A}}^{ *} {\rm {\bf A}}} \right)_{\,. \,i}
\left( {{\rm {\bf e}}_{.\,l} } \right)} \,\,\widetilde{d}_{lt}
\]
the $t$th component  of a row-vector ${\rm {\bf d}}^{{\rm {\bf
A}}}_{i\,.}= (d^{{\rm {\bf A}}}_{i1},...,d^{{\rm {\bf A}}}_{ip})$
for all $t=1,...,p$ and substitute it in (\ref{eq:x_ij_2}), we
obtain
\[x_{ij} ={\frac{{ \sum\limits_{t = 1}^{p}
 d^{{\rm {\bf A}}}_{it}
}{\sum\limits_{\alpha \in I_{r_{2},p} {\left\{ {j} \right\}}} {\rm
rdet}_{j} { \left( {({\rm {\bf B}}{\rm {\bf B}}^{ *} )_{j\,.\,}
({\rm {\bf e}}_{t.} )} \right)\,_{\alpha} ^{\alpha} } } }{{{{\rm
ddet} { {{\rm {\bf A}}} }} }{{\sum\limits_{\alpha \in I_{r_{2},p}}
{{\left| {\left( {{\rm {\bf B}}{\rm {\bf B}}^{ *} } \right){\kern
1pt} _{\alpha} ^{\alpha} } \right|}}} }}    }.
\]
Since again $\sum\limits_{t = 1}^{p}
 d^{{\rm {\bf A}}}_{it}{\rm {\bf e}}_{t.}={\rm {\bf
d}}^{{\rm {\bf A}}}_{i\,.}$, then it follows
(\ref{AXB_detA*A_d^A}), where ${\rm {\bf d}}^{{\rm {\bf
A}}}_{i\,.}$ is (\ref{eq:def_d^A}).

(iiii) The proof is similar to the proof of (iii). $\blacksquare$
\begin{corollary}(Theorem 3.3 in \cite{ky2})
Suppose
\begin{equation}\label{eq:two-sid}
{\rm {\bf A}}{\rm {\bf X}}{\rm {\bf B}} = {\rm {\bf C}}
\end{equation}
\noindent is a two-sided matrix equation, where ${\left\{ {{\rm
{\bf A}},{\rm {\bf B}},{\rm {\bf C}}} \right\}} \in {\rm M}(n,{\rm
{\mathbb{H}}} )$ are given, ${\rm {\bf X}} \in {\rm M}(n,{\rm
{\mathbb{H}}} )$ is unknown. If ${\rm ddet} {\rm {\bf A}} \ne 0$
and ${\rm ddet} {\rm {\bf B}} \ne 0$ , then (\ref{eq:two-sid}) has
a unique solution, and the solution is
\begin{equation}
\label{kyr10} x_{i\,j} = {\frac{{{\rm rdet} _{j} ({\rm {\bf
B}}{\rm {\bf B}}^{ *} )_{j.\,} \left( {{\rm {\bf c}}_{i\,.}^{{\rm
{\bf A}}}} \right)}}{{{\rm ddet}  {\rm {\bf A}}\cdot {\rm ddet}
{\rm {\bf B}}}}},
\end{equation}

\noindent or
\begin{equation}
\label{kyr11} x_{i\,j} = {\frac{{{\rm cdet} _{i} ({\rm {\bf A}}^{
*} {\rm {\bf A}})_{.\,i\,} \left( {{\rm {\bf c}}_{.j}^{{\rm {\bf
B}}}} \right)}}{{{\rm ddet} {\rm {\bf A}}\cdot {\rm ddet}  {\rm
{\bf B}}}}},
\end{equation}

\noindent where ${\rm {\bf c}}_{i\,.}^{{\rm {\bf A}}} : = \left(
{{\rm cdet} _{i} ({\rm {\bf A}}^{ *} {\rm {\bf A}})_{.\,i\,}
\left( {\tilde{{\rm {\bf c}}}_{.1}} \right),\ldots ,{\rm cdet}
_{i} ({\rm {\bf A}}^{ *} {\rm {\bf A}})_{.\,i\,} \left(
{\tilde{{\rm {\bf c}}}_{.n}} \right)} \right)$ is the row vector
and  ${\rm {\bf c}}_{.j}^{{\rm {\bf B}}} : = \left( {{\rm rdet}
_{j} ({\rm {\bf B}}{\rm {\bf B}}^{ *} )_{j.\,} \left( {\tilde{{\rm
{\bf c}}}_{1\,.}^{}} \right),\ldots ,{\rm rdet} _{j} ({\rm {\bf
B}}{\rm {\bf B}}^{ *} )_{j.\,} \left( {\tilde{{\rm {\bf
c}}}_{n\,.}^{}} \right)} \right)^{T}$ is the column vector and
$\tilde{{\rm {\bf c}}}_{i\,.}$, $\tilde{{\rm {\bf c}}}_{.\,j}$ are
the  ith row vector and the jth column vector of ${\rm {\bf
\widetilde{C}}}$, respectively, for all $i,j = 1,...,n $.
\end{corollary}
\begin{remark} In Eq.(\ref{eq:d^B}), the index $i$ in ${{\rm{cdet}} _{i} \left( {\left( {{\rm {\bf A}}^{
*} {\rm {\bf A}}} \right)_{\,.\,i} \left( {{{\rm {\bf
d}}}\,_{.\,j}^{{\rm {\bf B}}}} \right)} \right) _{\beta} ^{\beta}
}$ designates $i$th column of $\left( {\left( {{\rm {\bf A}}^{ *}
{\rm {\bf A}}} \right)_{\,.{\kern 1pt} \,i} \left( {{\rm {\bf
a}}_{.j}^{ *} } \right)} \right)$, but in the submatrix $\left(
{\left( {{\rm {\bf A}}^{ *} {\rm {\bf A}}} \right)_{\,.{\kern 1pt}
\,i} \left( {{\rm {\bf a}}_{.j}^{ *} }  \right)} \right){\kern
1pt} {\kern 1pt} _{\beta} ^{\beta} $ the entries of ${{\rm {\bf
a}}_{.j}^{ *} }$ may be placed in a column with the another index.
Similarly, we have for  $j$ in (\ref{eq:d^A}). In (\ref{eq:ls_AX})
and (\ref{eq:ls_XA}) we have equivalently.
\end{remark}
\section{An example}
In this section, we give an example to illustrate our results. Let
us consider the matrix equation
\begin{equation}\label{eq_ex:AXB=D}
 {\rm {\bf A}}{\rm {\bf X}}{\rm {\bf B}} = {\rm {\bf
D}},
\end{equation}
where
\[{\bf A}=\begin{pmatrix}
  1 & i & j \\
  -k & i & 1 \\
  k & j & -i \\
  j & -1 & i
\end{pmatrix},\,\, {\bf B}=\begin{pmatrix}
  i & 1 & j \\
 j & k & -i
\end{pmatrix},\,\, {\bf D}=\begin{pmatrix}
  1 & i & j \\
  k & 0 & i \\
  1 & j & 0 \\
  0 & k & i
\end{pmatrix}.\]
Then we have

\[{\bf A}^{*}{\bf A}=\begin{pmatrix}
  4 & 2i+2j & 2j+2k \\
  -2i-2j & 4 & -2i-2k \\
  -2j-2k & 2i+2k & 4
\end{pmatrix},\,\, {\bf B}{\bf B}^{*}=\begin{pmatrix}
  3 & -3k \\
  3k & 3
\end{pmatrix},\]
\[ {\rm {\bf \widetilde{D}}}= {\rm {\bf
A}}^\ast{\rm {\bf D}}{\rm {\bf B}}^{\ast}=\begin{pmatrix}
  2+2i+2j & -i+2j-2k \\
 1 -i-j & i-j-k \\
  1-2j & 2i-k
\end{pmatrix}.\]
Since ${\rm ddet}\,{\rm {\bf A}} =\det {\bf A}^{*}{\bf A}=0$ and
$\det{\bf A}^{*}{\bf A}=\begin{pmatrix}
  4 & 2i+2j  \\
  -2i-2j & 4
\end{pmatrix}=8\neq 0$, then  ${\rm rank}\,{\rm {\bf A}} = 2$.  Similarly ${\rm rank}\,{\rm {\bf
B}} = 1$.

So, we have the case (i) of Theorem \ref{theor:AXB=D}. We shall
find the minimum norm least squares solution ${\rm {\bf X}}_{LS}$
of (\ref{eq_ex:AXB=D}) by (\ref{eq:d^B}). We obtain \[
{{{\sum\limits_{\alpha \in I_{1,\,2}} {{\left| {\left( {{\rm {\bf
B}}{\rm {\bf B}}^{ *} } \right) {\kern 1pt} _{\alpha} ^{\alpha} }
\right|}}} }}=3+3=6,\]
\[ \begin{array}{c}
   {{{\sum\limits_{\beta \in J_{2,\,3}} {{\left| {\left( {{\rm {\bf
A}}^{ *}{\rm {\bf A}} } \right) {\kern 1pt} _{\beta} ^{\beta} }
\right|}}} }}=
  \det\begin{pmatrix}
  4 & 2i+2j \\
  -2i-2j & 4
\end{pmatrix}+\\\det\begin{pmatrix}
  4 & 2j+2k
 \\
 -2j-2k & 4
\end{pmatrix}+\det\begin{pmatrix}
  4 & -2i-2k \\
  2i+2k & 4
\end{pmatrix}=24.
\end{array}\]
By (\ref{eq:def_d^B_m}), we can get \[{\bf d}_{.1}^{{\bf
B}}=\begin{pmatrix}
  2+2i+2j \\
  1-i-j \\
  1-2j
\end{pmatrix},\,\,\,\,\,{\bf d}_{.2}^{{\bf  B}}=\begin{pmatrix}
  -i+2j+2k \\
  i-j-k \\
  2i-k
\end{pmatrix}.\]
Since
\[\left( {{\rm {\bf A}}^{ *} {\rm {\bf A}}} \right)_{\,.\,1}
\left( {{{\rm {\bf d}}}\,_{.\,1}^{{\rm {\bf B}}}} \right)=
\begin{pmatrix}
  2+2i+2j & 2i+2j & 2j+2k \\
  1-i-j & 4 & -2i-2k \\
  1-2j & 2i+2k & 4
\end{pmatrix},\] then finally we obtain
\[
  x_{11} = {\frac{{{\sum\limits_{\beta \in J_{2,\,3} {\left\{ {1}
\right\}}} { {\rm cdet} {\left(\left( {{\rm {\bf A}}^{ *} {\rm
{\bf A}}} \right)_{\,.\,1} \left( {{{\rm {\bf d}}}\,_{.\,1}^{{\rm
{\bf B}}}} \right)\right)\, _{\beta} ^{\beta}}  } }
}}{{{\sum\limits_{\beta \in J_{2,3}} {{\left| {\left( {{\rm {\bf
A}}^{ *} {\rm {\bf A}}} \right)_{\beta} ^{\beta} } \right|}}
\sum\limits_{\alpha \in I_{1,2}}{{\left| {\left( {{\rm {\bf
B}}{\rm {\bf B}}^{ *} } \right) _{\alpha} ^{\alpha} } \right|}}}
}}}=\]
 \[ \frac{{\rm cdet}_{1}\begin{pmatrix}
  2+2i+2j & 2i+2j \\
  1-i-j & 4
\end{pmatrix}+{\rm cdet}_{1}\begin{pmatrix}
  2+2i+2j & 2j+2k \\
  1-2j & 4
\end{pmatrix}}{144}=\]

\[\frac{4+5i+6j-k}{72}.
\]
Similarly,
\[
  x_{12} = \frac{{\rm cdet}_{1}\begin{pmatrix}
  -i+2j-2k & 2i+2j \\
  i-j-k & 4
\end{pmatrix}+{\rm cdet}_{1}\begin{pmatrix}
  -i+2j-2k & 2j+2k \\
  2i-k & 4
\end{pmatrix}}{144}=\]

\[\frac{-1-2i+5j-4k}{72},\]

    \[x_{21} =\frac{{\rm cdet}_{2}\begin{pmatrix}
  4 &  2+2i+2j \\
  -2i-2j & 1-i-j
\end{pmatrix}+{\rm cdet}_{1}\begin{pmatrix}
  1-i-j & -2i-2k \\
  1-2j & 4
\end{pmatrix}}{144}=\]

\[\frac{i-2j-k}{72},\]
\[
   x_{22} =\frac{{\rm cdet}_{2}\begin{pmatrix}
  4 & -i+2j-2k \\
  -2i-2j & i-j-k
\end{pmatrix}+{\rm cdet}_{1}\begin{pmatrix}
   i-j-k & -2i -2k\\
  2i-k & 4
\end{pmatrix}}{144}=\]

\[\frac{-2+2i+j-k}{72},\]

\[   x_{31} =\frac{{\rm cdet}_{2}\begin{pmatrix}
  4 & 2+2i+2j \\
 - 2j-2k & 1-2j
\end{pmatrix}+{\rm cdet}_{2}\begin{pmatrix}
   4 & 1-i-j\\
  2i+2k & 1-2j
\end{pmatrix}}{144}=\]

\[\frac{1-4i-3j}{72},\]
\[
    x_{32} =\frac{{\rm cdet}_{2}\begin{pmatrix}
  4 & -i+2j-2k \\
  -2j-2k &2i-k
\end{pmatrix}+{\rm cdet}_{1}\begin{pmatrix}
  4 &i-j-k \\
  2i+2k & 2i-k
\end{pmatrix}}{144}=\]

\[\frac{3i-3j-2k}{72}.
\]

Then
\[{\rm {\bf X}}_{LS}=\frac{1}{72}\left(
                                   \begin{array}{cc}
                                     4+5i+6j-k & -1-2i+5j-4k \\
                                     i-2j-k & -2+2i+j-k \\
                                     1-4i-3j & 3i-3j-2k \\
                                   \end{array}
                                 \right)\]
is the minimum norm least squares solution of (\ref{eq_ex:AXB=D}).

\end{document}